*Research Article - Preprint*

**Trust-region Filter Algorithms utilising Hessian Information for Grey-Box Optimisation**

Gul Hameed [1], Tao Chen [1], Antonio del Rio Chanona [2], Lorenz T Biegler [3], Michael Short [1*]

[1] School of Chemistry and Chemical Engineering, University of Surrey, Guildford GU2 7XH, UK
[2] Department of Chemical Engineering, Imperial College London, United Kingdom
[3] Department of Chemical Engineering, Carnegie Mellon University, Pittsburgh, PA 15213, USA

*Corresponding Author: m.short@surrey.ac.uk

**Abstract**

Optimising industrial processes often involves grey-box models that couple algebraic glass-box equations with black-box components lacking analytic derivatives. Such hybrid systems challenge derivative-based solvers. The classical trust-region filter (TRF) algorithm provides a robust framework but requires extensive parameter tuning and numerous black-box evaluations. This work introduces four Hessian-informed TRF variants (A1-A4) that use projected positive definite Hessians for automatic step scaling and minimal tuning, combined with both low-fidelity (linear, quadratic) and high-fidelity (Taylor series, Gaussian process) surrogates for local black-box approximation. Tested on 25 grey-box benchmarks and five engineering case studies (Himmelblau, liquid-liquid extraction, pressure vessel design, alkylation, and spring design), the new variants achieved up to an order-of-magnitude reduction in iterations and black-box evaluations, with reduced sensitivity to tuning parameters relative to the classical TRF algorithm. High-fidelity surrogates solved 92-100 % problems, compared to 72-84 % for the low-fidelity surrogates. Developed TRF methods also outperformed classical derivative-free optimisation solvers. The results show that new variants offer robust and scalable alternatives for grey-box process systems optimisation.

## 1. Introduction

Process systems optimisation has advanced considerably over the past few decades, enabling simultaneous optimisation and simulation of multiple process units. Early simulators adopted a sequential strategy, solving units in flowsheet order [1]. While easy to initialise and robust for small systems, this approach proved inefficient for large-scale optimisation, where unreliable derivative estimates often caused gradient-based solvers to fail or diverge. Modern simulators instead employ equation-oriented (EO) nonlinear models with automatic differentiation, providing reliable derivative information and allowing nonlinear programming (NLP) solvers to guarantee local optimality [1]. Widely used NLP solvers are based on active-set (CONOPT [2], MINOS [3], SNOPT [4]) and primal-dual interior point methods (KNITRO [5], LOQO [6], IPOPT [7]).

In practice, EO representations are not always available, and analytic derivative information may be inaccessible or prohibitively expensive. This situation arises when embedding high-fidelity simulations (e.g., computational fluid dynamics models from ANSYS Fluent) or proprietary modules (e.g., Aspen) within EO process flowsheets. For example, consider an EO optimisation model comprising multiple units, where some units are best represented by a computational fluid dynamics simulation. While a global surrogate could be constructed using machine learning, this approach typically requires excessive computational resources to achieve accuracy over the entire design space [8]. A more practical strategy is to treat these simulated units as external black-box calls during optimisation, constructing a local surrogate around the current iterate using a limited number of simulations, thereby avoiding the cost of a global surrogate. The resulting optimisation problem, combining a fully specified EO (glass-box) model with a black-box model lacking analytic derivatives, is termed grey-box. Here, EO models can be solved with deterministic NLP algorithms, while black-box modules often require pure derivative-free optimisation (DFO) methods such as direct search, Nelder-Mead, line-search, and trust-region strategies [9]. However, integrating deterministic and DFO approaches remains challenging because classical optimisation theory relies on analytic derivatives.

DFO methods struggle when constraints interact with black-box units, especially in recycle loops [10][11]. Each time decision variables change, the recycle loop must be fully converged, requiring many costly black-box evaluations per iteration, which makes DFO approaches computationally prohibitive for realistic process systems [10]. Large-scale and expensive grey-box models solved with pure DFO often fail to exploit available EO information, treating the entire problem as a black-box, leading to high computational costs and poor scalability [10]. Pure DFO also struggles with noisy or non-smooth objective functions and offers limited convergence guarantees without restrictive assumptions [12]. Despite advancements in modern DFO solvers such as COBYLA [13], BOBYQA [14], UOBYQA [15], NEWUOA [16], NOMAD [17], DAKOTA [18], and COBYQA [19], scalability to high-dimensional, highly non-linear, or constrained grey-box problems remains limited [12][20].

Among recent advancements in grey-box optimisation, the classical trust-region filter (TRF) method has emerged as a promising approach [1][10][21]. TRF combines surrogate modelling, trust-region methods, filter strategy, DFO principles, and implicit handling of constraints (via nonlinear solver) [9][22]. Originating from the SQP-filter technique [23] and its extensions to reduced-order models (RMs) [24][25] and inexact Jacobians [26], TRF was first applied to partial-differential-constrained optimisation of pressure swing adsorption (PSA) $CO_2$ capture using proper orthogonal decomposition (POD) RMs [25]. The method was later extended to grey-box optimisation with proven global convergence under mild assumptions [1][10]. Applications include optimisation of the Williams-Otto, ammonia synthesis, air-fired power plant, and $CO_2$ capture process flowsheets. TRF strategies have also been adapted to handle high-fidelity [27] or rigorous models [28] (with available sensitivity information) using first-order Taylor series (TS) approximations within nonlinear EO models. More recently, probabilistic Gaussian Process (GP) surrogates have been incorporated into TRF for shape optimisation of biocatalytic micro-reactors, where convergence relied solely on local surrogate accuracy [29]. Nonetheless, review of TRF methods [29] highlights challenges, including high iteration counts, many black-box calls, and sensitivity to tuning parameters [30], making the existing TRF algorithms expensive for grey-box problems involving costly NLP models or expensive black-box simulations.

Beyond the development of TRF strategies, several significant data-driven and surrogate-based frameworks have shaped grey-box optimisation in process systems engineering. Pistikopoulos and Ierapetritou were among the first to propose decomposition-based strategies coupling feasibility subproblems (defining the feasible region) and a master optimality problem (balancing feasibility and optimality) for process optimisation under model uncertainty [31]. Later hybrid simulation-optimisation approaches were developed, linking process simulators with optimisation models [32][33]. Wang and Ierapetritou proposed a black-box simulation-optimisation method using stochastic Kriging surrogates with sequential sampling to approximate objectives and constraints [34]. These early frameworks established the foundation for hybrid modelling but lacked convergence guarantees and primarily targeted stochastic or dynamic simulation/black-box problems rather than large-scale static nonlinear grey-box flowsheets. ARGONAUT (Algorithms for Global Optimisation of Constrained Grey-box Computational Problems) combined variable selection, bounds tightening, and constrained sampling to construct surrogate representations of unknown equations, which were then globally optimised [35]**.** For constrained cases, ARGONAUT employed a two-stage sampling strategy, first generating a large Latin Hypercube design within bound-constrained space, then filtering infeasible samples and optimally selecting a reduced subset based on probabilistic input-output distance criteria. While this constrained sampling procedure improved surrogate quality, it increased computational cost for large systems. The latter p-ARGONAUT implementation parallelised these stages [36], but required high-performance computing resources. ALAMO (Automated Learning of Algebraic Models for Optimisation) introduced an adaptive regression framework to construct algebraic surrogates from simulation data, with constraint-handling capabilities and systematic improvement using derivative-free optimisation [37][38]. However, its ability to capture strong nonlinearities and accurate model behaviour was limited by the predefined basis-function library. Bajaj and Hasan developed a deterministic global derivative-free algorithm exploiting bounded-Hessian assumptions, bounding the diagonal elements of the Hessian using physical knowledge, to guarantee convergence for smooth black-box functions [39]. While rigorous, this approach applied only to problems where valid lower bounds could be computed, required numerous black-box evaluations, and treated the entire problem as a black box, excluding glass-box integration. Trust-region methods using radial basis function (RBF) [40][41] or Kriging surrogates [42], and an inexact restoration approach relying on polynomial models [43], further advanced DFO algorithms but rarely included second-order information or NLP globalisation strategies such as a filter. Recent reviews highlight rapid adoption of machine-learning surrogates such as GPs and neural networks [44][33], yet note that few frameworks integrate explicit Hessian information or simultaneously address curvature, feasibility, and optimality within a unified framework. It has also been emphasised that different optimisation problems may require different surrogate forms [33][45][46], whereas most existing algorithms are tied to specific surrogate types.

Collectively, these developments represent major progress yet reveal persistent gaps. Global surrogate frameworks rely on extensive sampling and retraining; regression-based models neglect curvature information; bounded-Hessian methods require many black-box evaluations and apply only to fully black-box problems; many approaches lack surrogate flexibility or adaptability under convergence guarantees; and machine-learning surrogates often sacrifice theoretical rigour. This work aims to address these limitations by extending the classical TRF framework (Algorithm A0) [10] through four new constrained grey-box optimisation algorithms (A1-A4). These new variants integrate local Hessian projections to guide the optimisation, reducing external evaluations, accelerating convergence and minimising tuning parameter sensitivity [27][47]. The new variants also allow a unified surrogate modelling layer that enables switching between low-fidelity (linear, quadratic, simplified quadratic) and high-fidelity (Taylor series, GP) surrogate models depending upon the problem at hand. In summary, we present curvature-aware, surrogate-flexible TRF algorithms that reduce external evaluations, improve robustness, and preserve the global convergence guarantees of the classical method [21].

The remainder of this paper is organised as follows. Section 2 presents the mathematical formulation of the grey-box optimisation problem and the components of the classical TRF algorithm. Section 3 describes the proposed surrogate modelling strategies and Hessian-informed variants. Section 4 briefly mentions implementation details. Section 5 provides numerical and engineering case studies demonstrating performance improvements over the classical TRF and pure DFO methods, and discusses the broader implications of these findings. Finally, Section 6 concludes the paper with key insights and future research directions.

## 2. Algorithmic Concepts

The trust-region framework is widely used for NLPs by approximating objectives and constraints [23]. In the grey-box TRF method, the trust region controls the approximation error of black-box surrogates [48].

We consider a grey-box problem

$$\min_{w,z} f\big(z, w, d(w)\big) \ s.t. \ \ h\big(z, w, d(w)\big) = 0 \ and \ g(z, w, d(w)) \leq 0, \tag{1}$$

with $w \in \mathbb{R}^m$ (black-box inputs/variables), $z \in \mathbb{R}^n$ (remaining variables), and a black-box map $d(w) : \mathbb{R}^m \to \mathbb{R}^p$. Assume $f$, $h$, $g$ and $d$ are twice continuously differentiable ($C^2$) functions.

To decouple the glass-box and black-box terms, introduce an auxiliary variable $y \in \mathbb{R}^p$:

$$\min_{x} f(x) \ s.t. \ \ h(x) = 0, g(x) \leq 0 \ and \ y = d(w). \tag{2}$$

In all formulations, $x = (w, y, z) \in \mathbb{R}^{m+n+p}$ denotes the complete decision-variable vector.

At iteration $k$, a local surrogate model $s^{(k)}(w)$ of $d(w)$ is constructed within the closed ball $D\big(w^{(k)}, \Delta^{(k)}\big) = \big\{w \in \mathbb{R}^m \colon \big\|w - w^{(k)}\big\| \leq \Delta^{(k)}\big\}$ centered at $w^{(k)}$ with trust-region radius $\Delta^{(k)}$. To guarantee global convergence, surrogates $s^{(k)}(w)$ must be $\kappa - fully\ linear$, satisfying:

$$\big\|\nabla s^{(k)}(w) - \nabla d(w)\big\| \leq \kappa_g \Delta^{(k)} \ and \ \big\|s^{(k)}(w) - d(w)\big\| \leq \kappa_f {\Delta^{(k)}}^2, \tag{3}$$

for all $w$: $\big\|w - w^{(k)}\big\|$. As $\Delta^{(k)} \to 0$, the *κ-fully linear* property ensures $s^{(k)}(w)$ converges to the true model $d(w)$ for some $\kappa_g, \kappa_f > 0$ [9]. Even if $\nabla d(w)$ is unavailable; its existence is assumed for theoretical completeness. Surrogates are reconstructed if *κ-fully linear* property (3) is violated.

The surrogate and trust-region constraint results in a trust-region subproblem (TRSP):

$$\min_{x} f(x) \ s.t. \ h(x) = 0, g(x) \leq 0, y = s^{(k)}(w) \ and \ \big\|x - x^{(k)}\big\| \leq \Delta^{(k)}. \tag{4}$$

The TRF method iteratively solves the TRSP to generate a sequence $\{x^{(k)}\}$ converging to a first-order Karush-Kuhn-Tucker (KKT) point of the original grey-box problem (1), ensuring feasibility for both glass-box and black-box constraints, while requiring minimal black-box evaluations.

A TRF iteration comprises the following algorithmic components. Criticality check (section 2.1) and *κ-fully linear* property, ensure surrogate accuracy and check termination. Compatibility check (section 2.2) initialises and ensures the feasibility of the TRSP. The restoration phase (section 2.6) is invoked when the compatibility check fails; otherwise, the algorithm proceeds to the TRSP solution (section 2.3). Filter and trust-region update mechanisms (section 2.4) govern how the TRSP solution is used to advance to the next iteration.

### 2.1 Criticality Check

The criticality check, trust-region size and surrogate accuracy determine whether the algorithm has reached a KKT point of the grey-box problem. The criticality measure, $\chi^{(k)}$, is obtained by solving a linearised TRSP (i.e., the criticality problem) at the current iterate:

$$\chi^{(k)} = \left|\min_{v} \nabla f\left(x^{(k)}\right)^T v\right| \quad (5)$$

$$s.t. \nabla h\left(x^{(k)}\right)^T v = 0, g\left(x^{(k)}\right) + \nabla g\left(x^{(k)}\right)^T v \leq 0, v_y - \nabla s^{(k)}\left(w^{(k)}\right)^T v_w = 0 \; and \; \|v\| \leq 1.$$

Here $v = \left(v_w, v_y, v_z\right)$ is a unit trust-region step, analogous to $x = (w, y, z)$. If a polyhedral norm ($l_\infty$ or $l_1$) is used in $\|v\| \leq 1$, (5) becomes a linear program.

A small $\chi^{(k)}$ value indicates proximity to stationarity. If $\chi^{(k)}$ falls below a threshold relative to trust-region size (i.e., $\chi^{(k)} < \xi\Delta^{(k)}$ for $\xi > 0$), $\Delta^{(k)}$ is reduced ($\Delta^{(k+1)} = \omega\Delta^{(k)}$ for $0 < w < 1$). By iteratively shrinking $\Delta^{(k)}$, the algorithm enforces $\chi^{(k)} \rightarrow 0$ and guarantees accurate derivatives via $\kappa - fully\ linear$ property (3), thereby certifying optimality of the original grey-box problem [9].

### 2.2 Compatibility Check

The trust-region size constraint and the surrogate may render the candidate TRSP infeasible. When the feasible region lies close $x^{(k)}$, the subproblem remains a good local approximation of the original problem.

Feasibility is assessed by solving the compatibility problem:

$$\min_{w,z} \left\|y - s^{(k)}(w)\right\| \; s.t. \;\; h(x) = 0, g(x) \leq 0 \; and \; \left\|x - x^{(k)}\right\| \leq \kappa_\Delta \Delta^{(k)} min\left[1, \kappa_\mu {\Delta^{(k)}}^{\mu}\right], \quad (6)$$

where $y$ is the current fixed surrogate output, and $\mu \in (0,1)$, $\kappa_\mu > 0$, and $\kappa_\Delta \in (0,1)$ are constants. This problem is always feasible at $x = x^{(k)}$, since $h\left(x^{(k)}\right) = 0$ and $g(x^{(k)}) \leq 0$.

Let $\iota^{(k)} = \left\|y - s^{(k)}(w)\right\|$. If $\iota^{(k)} = 0$ (or $\iota^{(k)} \leq \epsilon_{comp}$), candidate TRSP is compatible; otherwise, the current iterate is added to the filter $\mathcal{F}^{(k)}$ (section 2.4) and a restoration phase (section 2.6) is initiated to generate a TRSP-compatible surrogate and the trust-region size.

The compatible step $\bar{r}^{(k)}$ is:

$$\bar{r}^{(k)} = x - x^{(k)} = \left(\bar{r}_w^{(k)}, \bar{r}_y^{(k)}, \bar{r}_z^{(k)}\right), \quad (7)$$

analogous to the normal step in composite-step NLP methods. If the compatibility holds, $\left\|\bar{r}^{(k)}\right\| \leq \kappa_{usc}\theta^{(k)}$ for finite $\kappa_{usc}$ (independent of $k$), where $\theta^{(k)}$ is the infeasibility measure as defined by equation (10).

### 2.3 Trust-Region Subproblem

The TRSP is initialised at $x^{(k)} + \bar{r}^{(k)}$ and solved to obtain a trial step. The resulting step is:

$$r^{(k)} = x_r^{(k)} - x^{(k)}, \tag{8}$$

where $x_r^{(k)}$ is the minimiser of the TRSP.

The step $r^{(k)}$ must satisfy the fraction-of-Cauchy-decrease condition for the convergence [1]:

$$f\left(x^{(k)} + \bar{r}^{(k)}\right) - f\left(x^{(k)} + r^{(k)}\right) \geq \kappa_{tmd}\chi^{(k)} min\left[\frac{\chi^{(k)}}{\beta^{(k)}}, \Delta^{(k)}\right], \tag{9}$$

where $\kappa_{tmd}$ is a constant and $\left\{\beta^{(k)}\right\}$ is a bounded sequence with $\beta^{(k)} > 1$.

### 2.4 Filter Mechanism and Trust-Region Update

The filter mechanism, originally developed for the filter-SQP method [23], has been adapted in TRF-based grey-box optimisation [1]. Unlike penalty-based methods for NLPs (e.g., $P(x, \mu) = f(x) + \zeta h'(c)$, where $P(x, \mu)$ is a penalty function, $h'(c)$ is a measure of constraint violation and $\zeta > 0$ is a penalty parameter), which require careful parameter tuning, the filter approach reformulates the optimisation as a bi-objective problem [23]. It simultaneously minimises the objective and either the measure of constraint violation (for NLPs) or the infeasibility (for grey-box problems).

In TRF methods, the infeasibility measure is defined as:

$$\theta(x) = \|s(w) - d(w)\|. \tag{10}$$

At certain iterations $k$, pairs $(f^{(k)}, \theta^{(k)})$ are added to the filter set

$$\mathcal{F}^{(k)} = \left\{\left(f^{(j)}, \theta^{(j)}\right): j < k, j \in \mathcal{Z}\right\}, \tag{11}$$

where $\mathcal{Z} \subset \mathbb{N}$ indexes filter iterations. This builds a Pareto front of non-dominated $\left(f^{(k)}, \theta^{(k)}\right)$ pairs (i.e., those that improve both objective and infeasibility measure compared to existing filter entries), see Figure 1.

A trial step $r^{(k)}$ is accepted if the filter condition holds:

$$\theta\left(x^{(k)} + r^{(k)}\right) \leq (1 - \gamma_\theta)\theta^{(j)} \quad or \quad f\left(x^{(k)} + r^{(k)}\right) \leq f^{(j)} - \gamma_f\theta^{(j)}, \tag{12}$$

with $\gamma_\theta \in (0,1)$ and $\gamma_f \in (0,1)$.

**Step Classification**

The accepted steps are categorised using the switching condition:

$$f(x^{(k)}) - f(x^{(k)} + s^{(k)}) \geq \kappa_\theta \theta(x^{(k)})^{\gamma_s} and\ \theta^{(k)} \leq \theta_{min}, \tag{13}$$

where $\gamma_s > \frac{1}{1+\mu}$ and $\kappa_\theta \in (0,1)$.

- If (13) holds, the step is $f - type$, meaning that the surrogate error is acceptable and the trust-region radius expands as $\Delta^{(k+1)} := \max[\gamma_e \| r^{(k)} \|, \Delta^{(k)}]$, with $\gamma_e > 1$.
- Otherwise, the step is $\theta - type$, the pair $(f^{(k)}, \theta^{(k)})$ is added to the filter set, and the trust-region size is updated via the ratio test (14):

$$\Delta^{(k+1)} = \begin{cases} \gamma_c \| r^{(k)} \| & if\ \rho^{(k)} < \eta_1, \\ \Delta^{(k)} & if\ \eta_1 \leq \rho^{(k)} < \eta_2, \\ \max[\gamma_e \| r^{(k)} \|, \Delta^{(k)}] & if\ \rho^{(k)} \geq \eta_2, \end{cases} \tag{14}$$

where $0 < \gamma_c < 1 < \gamma_e$ and $0 < \eta_1 \leq \eta_2 < 1$. According to convergence theory, there is flexibility in updating $\Delta^{(k)}$ for $\theta - type$ steps provided the radius remains bounded [24].

The ratio of the actual reduction of the infeasibility measure to its predicted reduction, i.e., $\rho^{(k)}$, is given by:

$$\rho^{(k)} = \frac{\theta(x^{(k)}) - \theta(x^{(k)} + r^{(k)}) + \epsilon_\theta}{\max(\| d(w^{(k)}) - s^{(k)}(w^{(k)}) \|, \epsilon_\theta)}, \tag{15}$$

where $\epsilon_\theta$ is a small positive tolerance, preventing the decrease in trust-region size when $\theta(x^{(k)})$ and $\theta(x^{(k)} + r^{(k)})$ are very small.

For convergence, the trust-region radius must not decrease when the step is $f - type$, but must decrease when a step is rejected [10]. Beyond these fundamental requirements, the trust-region update method can be modified using sensitivity analysis or problem-specific heuristics.

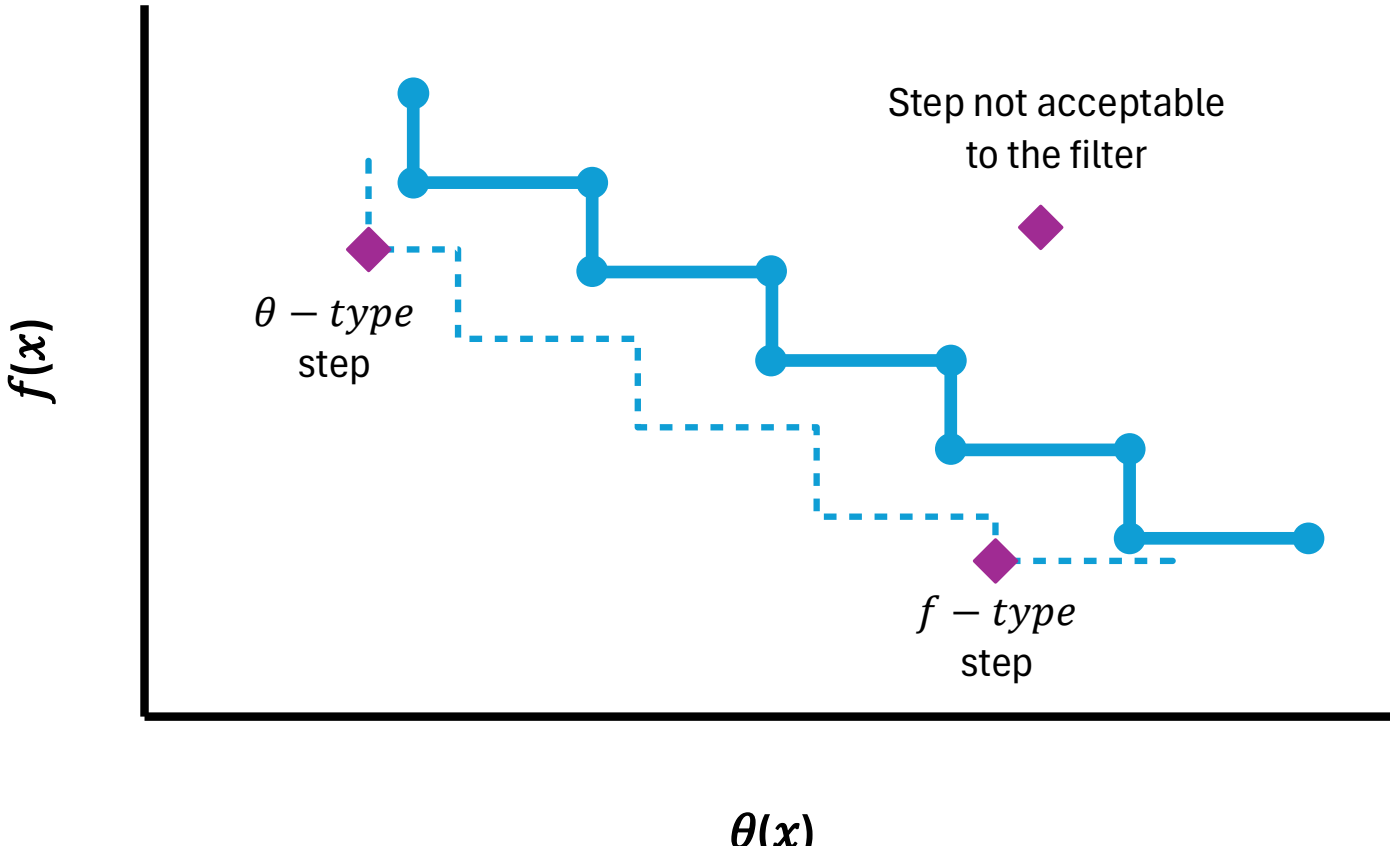

***Figure 1.*** *Illustration of the filter mechanism used in the TRF framework. Solid blue dots represent previously accepted* $(\boldsymbol{f},\boldsymbol{\theta})$ *pairs in the filter, dashed blue lines indicate the acceptance boundary defined by filter condition (12), and purple diamonds show trial points tested against this criterion. The filter accepts a step if it improves either objective* ($\boldsymbol{f - type}$) *or feasibility* ($\boldsymbol{\theta - type}$)*, thereby guiding progress toward both feasibility and optimality without penalty parameters.*

### 2.5 Sampling-region

In trust-region-based DFO, the trust-region serves a dual role: limiting step size and defining the region where the surrogate is accurate. While this worked well in pure black-box settings (where model accuracy and step control are aligned), it was noted that these roles should be decoupled. In UOBYQA [15], a separate sampling-region was introduced to construct $\kappa - fully\ linear$ surrogate models around black-box variables, while the trust-region continued to govern globalisation. This separation is important in grey-box problems, where the traditional trust-region constraint applies to both glass-box and black-box variables, restricting progress in glass-box variables that do not directly interact with the black-box model [1].

At iteration $k$, the sampling-region is defined as a space $S(w^{(k)}, \sigma^{(k)})$, where $\sigma^{(k)}$ is sampling-region radius ($\sigma^{(k)} \leq \Delta^{(k)}$), and $w^{(k)}$ is the centre (i.e., the recently accepted iterate). This ensures surrogate models remain valid even if the TRSP step is rejected. Unlike traditional methods, where $\Delta^{(k)} \rightarrow 0$ for convergence, here only $\sigma^{(k)}$ contracts, allowing glass-box variables to take larger steps even near optimality.

During criticality updates, shrinking the sampling-region is preferred over shrinking the trust-region. The standard criticality test is updated to:

$$\chi^{(k)} < \xi \sigma^{(k)}. \tag{16}$$

If (16) holds, $\sigma^{(k)}$ is updated as:

$$\sigma^{(k)} = \max\left(\min\left(\sigma^{(k-1)}, \chi^{(k)}/\xi\right), \Delta_{min}\right). \tag{17}$$

The criticality check (16) and $\kappa - fully\ linear$ property (with $\sigma$ substituting for $\Delta$ in (3)) preserve the first-order optimality.

### 2.6 Restoration Phase

If the compatibility check fails, a restoration phase is invoked to re-establish the TRSP feasibility. Starting from the current iterate $x^{(k)}$ and radius $\Delta^{(k)}$, the compatibility problem (6) is solved iteratively while updating $\Delta^{(k)}$ via the reduction ratio test (14). This process continues until a compatible set ($x^{(k+1)}$, $\Delta^{(k+1)}$, $s^{(k+1)}$) is obtained. The restored iterate is guaranteed to be feasible and always satisfies the filter condition (12). Restoration ensures that subsequent TRSPs are solvable and that convergence is maintained.

### 2.7 Existing Algorithms

Two TRF algorithms have previously been proposed by Eason and Biegler [10], differing in their treatment of surrogate accuracy and termination [10]. The first algorithm [1] adjusts the search step and enforces surrogate accuracy with a single trust-region. The second algorithm [10] employs surrogate accuracy within the sampling-region and introduces explicit termination conditions based on TRF convergence proofs.

Both algorithms were tested on the COPS test set [49] using linear (18) and quadratic (19) surrogate models:

$$s_{linear} = b_0 + \sum_{i=1}^{n} b_i w_i, \tag{18}$$

$$s_{quadratic} = b_0 + \sum_{i=1}^{n} b_i w_i + \sum_{i=1}^{n} b_i' w_i^2 + \sum_{1 \le i < j \le n} m_{ij} w_i w_j, \tag{19}$$

where $w_i$ are surrogate inputs, $b_0$ is a bias term (i.e., true black-box value at the centre of the trust-region) ensuring the *κ-fully linear* property, and $b_i$, $b_i'$ and $m_{ij}$ are polynomial coefficients obtained from black-box evaluations. Here, $n$ denotes surrogate dimension. The second algorithm [10] solved all problems in the COPS set, whereas the first [1] was only able to solve half. Therefore, the second algorithm (Figure 2) serves as the base classical TRF algorithm (A0) in this work.

**Base Algorithm (A0):**

1. **Initialisation:** Choose tuning parameters $0 < \gamma_c < 1 < \gamma_e$, $0 < \eta_1 \le \eta_2 < 1$, $\gamma_\theta \in (0,1)$, $\gamma_f \in (0,1)$, $\kappa_\theta \in (0,1)$, $\mu \in (0,1)$, $\gamma_s > 1/(1+\mu)$, $\kappa_\Delta \in (0,1)$, $\kappa_\mu > 1$, $\xi > 0$, $\Psi \in (0,1)$. Specify tolerances: $\epsilon_\theta$, $\epsilon_\chi$, $\epsilon_{comp}$, and $\epsilon_\Delta \ge \Delta_{min}$. Define bounds $\Delta_{min} > 0$, $\theta_{min} > 0$ and $\theta_{max} > 0$.

   Set $x^{(0)}$, $\Delta^{(0)} > 0$, $\sigma^{(0)} > 0$. Evaluate $d(w^{(0)})$, $\theta(x^{(0)})$. Initialise $\mathcal{F}^{(0)} = \emptyset$ and $k = 0$.
2. **Surrogate Construction:** Build a surrogate model $s^{(k)}(w^{(k)}w)$ that is *κ-fully linear* within a sampling-region of radius $\sigma^{(k)}$ centred at $x^{(k)}$.
3. **Termination and criticality check:** Calculate $\chi^{(k)}$ using the criticality problem (5).
   a. If $\theta^{(k)} \le \epsilon_\theta$, $\chi^{(k)} \le \epsilon_\chi$, and $\sigma^{(k)} \le \epsilon_\Delta$, STOP (first-order critical point).
   b. If $\theta^{(k)} \le \epsilon_\theta$, $\theta^{(k-1)} \le \epsilon_\theta$, $\Delta^{(k)} \le \Delta_{min}$, and $\Delta^{(k-1)} \le \Delta_{min}$, STOP (feasible point).
   c. Else if the criticality test (16) holds, perform criticality update (17).
4. **Compatibility check:** Solve the compatibility problem (6). If $\|y - s^{(k)}(w)\| \le \epsilon_{comp}$, go to the next step. Otherwise, add $(f^{(k)}, \theta^{(k)})$ into $\mathcal{F}^{(k+1)}$ and call restoration phase iteration, then set $k = k + 1$ and go to step 2. The restoration loop (steps 2-4) continues until compatible $x^{(k+1)}$, $\Delta^{(k+1)} > 0$, $\sigma^{(k+1)} > 0$, and $s^{(k+1)}(w)$ are found.
5. **Trust-region subproblem:** Initialise the TRSP (4) and solve it to compute the step $r^{(k)}$.
6. **Filter acceptance:** Evaluate $\theta^{(k)}$ and $f^{(k)}$.
   a. If filter condition (12) is satisfied, continue to Step 7.
   b. Else, set $x^{(k+1)} = x^{(k)}$, $\Delta^{(k+1)} = \gamma_c \|r^{(k)}\|$, $\theta^{(k+1)} = \theta^{(k)}$, $\sigma^{(k+1)} = \min\{\sigma^{(k)}, \Psi\Delta^{(k)}\}$, set $k = k + 1$ and go to step 2.
7. **Switching:** If the switching condition (13) holds, go to step 8. Else, go to step 9.
8. ***f* − *type* step:** Update $x^{(k+1)} = x^{(k)} + r^{(k)}$, $\Delta^{(k+1)} = \max\{\gamma_e \|r^{(k)}\|, \Delta^{(k)}\}$, $\sigma^{(k+1)} = \sigma^{(k)}$, $\theta^{(k+1)} = \theta(x^{(k)} + r^{(k)})$, set $k = k + 1$ and go to step 2.

9. $\boldsymbol{\theta - type}$ **step:** Add $\left(f^{(k)}, \theta^{(k)}\right)$ into $\mathcal{F}^{(k+1)}$. Update $x^{(k+1)} = x^{(k)} + r^{(k)}$, $\sigma^{(k+1)} = \min[\sigma^{(k)}, \Psi\Delta^{(k)}]$, adjust trust-region size using ratio test rule (14)-(15), set $k = k + 1$ and go to step 2.

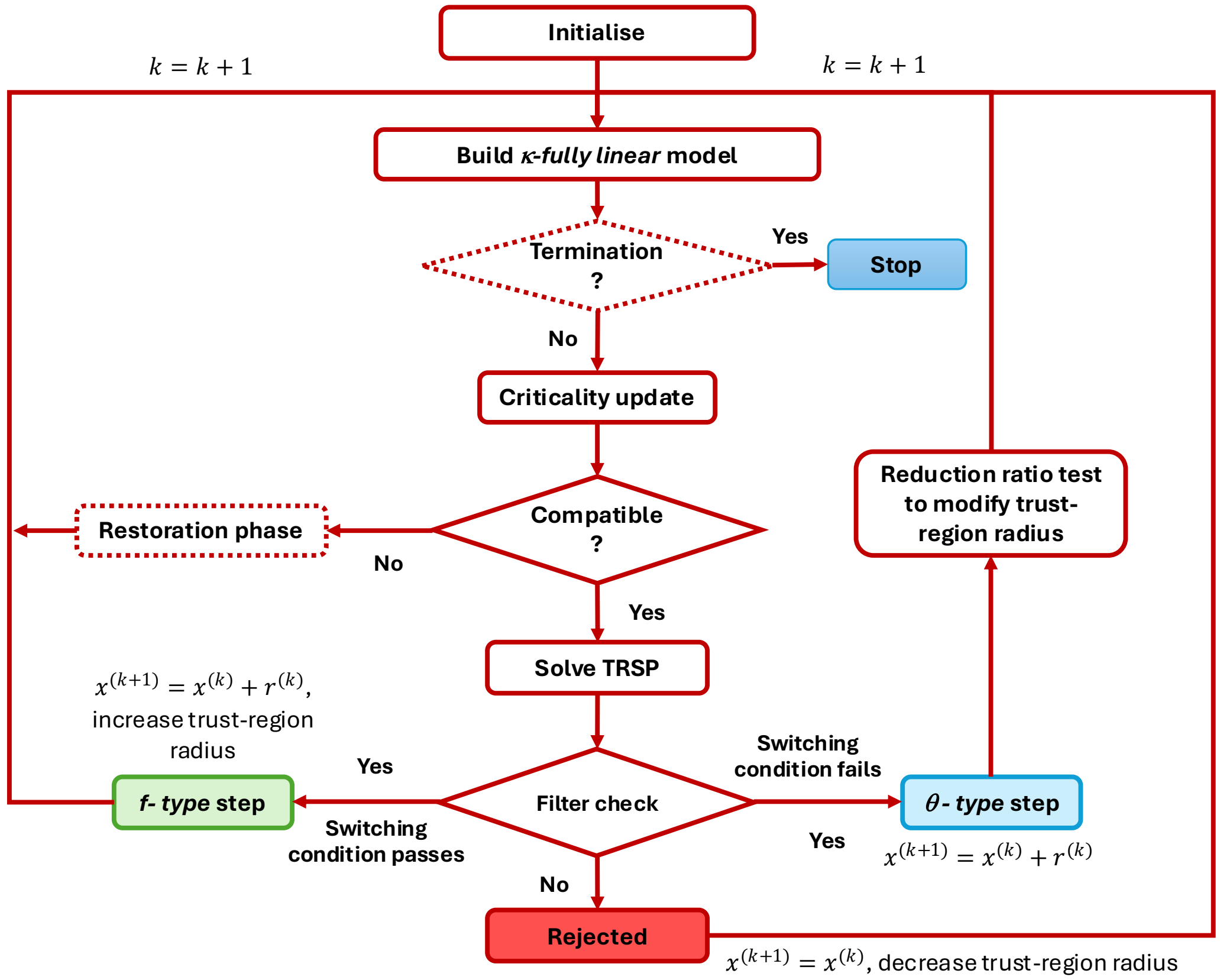

***Figure 2.*** *Base algorithm (A0) of the classical TRF method. The flowchart depicts the sequence of operations, including surrogate construction, criticality and compatibility checks, TRSP solution, and filter acceptance. The iterative loops involving restoration,* $f - type$*, and* $\theta - type$ *steps ensure global convergence to a first-order critical point.*

## 3. Proposed Improved Algorithms

Building on the trust-region theory [22], DFO [9], TRF developments for grey-box optimisation [27][29], and the generalised trust-region Newton method for highly non-convex hyperelastic simulations [50][51], we propose several enhancements to the base TRF algorithm (A0) to improve robustness, scalability, and efficiency. Improvements include surrogate modelling, incorporation of local Hessian information, and refinement of termination criteria.

### 3.1 Surrogate Modelling

The computational expense of grey-box optimisation is primarily due to constructing a $\kappa - fully\ linear$ surrogate. Since TRF methods require multiple black-box evaluations per iteration, efficient surrogate modelling is crucial for reducing costly external evaluations.

Polynomial surrogates are commonly used in TRF frameworks and inherently satisfy the $\kappa - fully\ linear$ property [10]. Quadratic polynomials (21) capture variable interactions but require $(n + 1)(n + 2)/2$ sample points, where $n$ is the dimension of the true model (i.e., the number of black-box inputs). In contrast, a simplified quadratic (20), which neglects interaction terms, requires only $2n + 1$ sample points [29].

$$s_{simp_quadratic} = b_0 + \sum_{i=1}^{n} b_i w_i + \sum_{i=1}^{n} b_i' w_i^2 \tag{20}$$

If both the function value $d\left(w^{(k)}\right)$ and gradient $\nabla d\left(w^{(k)}\right)$ at the current iterate are available as an external call, and a global basis function $s_{gb}(w)$ (e.g., polynomial, Kriging, machine learning model) can be constructed, then a first-order TS approximation $s_{ts}$ can be formed:

$$\begin{aligned} s_{ts} = s_{gb}(w) + \left(d\left(w^{(k)}\right) - s_{gb}\left(w^{(k)}\right)\right) \\ + \left(\nabla d\left(w^{(k)}\right) - \nabla s_{gb}\left(w^{(k)}\right)\right)^T \left(w - w^{(k)}\right). \end{aligned} \tag{21}$$

Note that $s_{gb}(w)$ can even be set to zero in (21).

While polynomial and TS surrogates can provide efficient $\kappa - fully\ linear$ models, they have limitations. They rely on local sample points and require well-poised sets for stable fitting. TS surrogates also require derivatives, which may not always be available. These limitations reduce flexibility when data are sparse or unevenly distributed. GP surrogates address these issues by providing a probabilistic, non-parametric framework. GPs naturally quantify uncertainty and weight sample points adaptively via a kernel, making them effective for modelling complex nonlinear response surfaces.

Given a training set $\mathcal{D}$ with training points $\boldsymbol{W}$ and corresponding observations $\boldsymbol{y}$ (22):

$$\mathcal{D} = \{\boldsymbol{W}, \boldsymbol{y}\}, \boldsymbol{W} = [\boldsymbol{w_1}, \dots, \boldsymbol{w_n}]^T \in \mathbb{R}^{n \times d}, \boldsymbol{y} = [y_1, \dots, y_n]^T \in \mathbb{R}^n, \tag{22}$$

we place a GP prior over the latent function in the functional space expressed as:

$$f(.) \sim \mathcal{GP}(m(.), k(.,.)), \tag{23}$$

where $m(.)$ is the mean function and $k(.,.)$ is the kernel function. The underlying kernel function plays a crucial role in defining the smoothness and structure of the surrogate model. Common choices such as the squared exponential kernel and Matérn kernels balance smoothness and adaptability [52].

Evaluating the GP prior (23) at the training inputs $\{w_1, \dots, w_n\}$ yields a joint Gaussian distribution over the corresponding function values $\{y_1, \dots, y_n\}$ [53]:

$$\boldsymbol{y} = [f(\boldsymbol{w_1}), \dots, f(\boldsymbol{w_n})]^T \sim \mathcal{N}(\boldsymbol{m}, \boldsymbol{K}), \tag{24}$$

with mean vector $\boldsymbol{m} = m(\boldsymbol{W}) = [m(\boldsymbol{w_1}), \dots, m(\boldsymbol{w_n})]$ and the covariance matrix $\boldsymbol{K_{nn}} = k(\boldsymbol{W}, \boldsymbol{W}) \in \mathbb{R}^{n \times n}$ with entries $k(w_i, w_j)$.

Conditioning on an unknown test input $\boldsymbol{w}^* \in \mathbb{R}^d$, we evaluate the prior joint distribution between the known training outputs ($\boldsymbol{y}$) and the unknown output ($f^* = f(\boldsymbol{w}^*)$) [54]:

$$\begin{bmatrix} \boldsymbol{y} \\ f^* \end{bmatrix} \backsim \mathcal{N}\left( \begin{bmatrix} m(\boldsymbol{W}) \\ m(\boldsymbol{w}^*) \end{bmatrix}, \begin{bmatrix} \boldsymbol{K_{nn}} & \boldsymbol{K_{n*}} \\ \boldsymbol{K_{n*}}^T & \boldsymbol{K_{**}} \end{bmatrix} \right), \tag{25}$$

where $\boldsymbol{K_{n*}} = k(\boldsymbol{w}^*, \boldsymbol{W})$ and $\boldsymbol{K_{**}} = k(\boldsymbol{w}^*, \boldsymbol{w}^*)$ is prior variance at the unknown test point.

Using standard results from conditioning Gaussian distributions, the posterior distribution over $f(\boldsymbol{w}^*)|\mathcal{D}$ is also Gaussian [29]:

$$f(\boldsymbol{w}^*)|\mathcal{D} \sim \mathcal{N}(m', k'), \tag{26}$$

with

$$m' = m(\boldsymbol{w}^*) + \boldsymbol{K_{n*}}^T \boldsymbol{K_{nn}}^{-\mathbf{1}}(\boldsymbol{y} - \boldsymbol{m}), \;\; k' = \boldsymbol{K_{**}} - \boldsymbol{K_{n*}}^T \boldsymbol{K_{nn}}^{-1} \boldsymbol{K_{n*}}. \tag{27}$$

Often $m(.) = 0$ is assumed for mean-centred data, though prior knowledge may be encoded through a parametric or data-driven mean function.

By assigning weights based on Euclidean distance, GP regression ensures that nearby observations contribute more to the model while filtering out distant/irrelevant points [29]. This enables a more efficient use of black-box evaluations, reducing external function calls. This may improve both the efficiency and robustness of grey-box optimisation, particularly when black-box evaluations are highly non-linear [29].

Finally, we propose a hybrid surrogate $s_{hybrid}$, which augments a TS surrogate with a GP residual model. Specifically, a Taylor expansion (zero- and first-order terms, as in Eq. 21) is combined with a GP trained on the residuals $s_{GP,res}$ (i.e., the difference between true function values and the Taylor approximation):

$$s_{hybrid} = s_{ts}\left(w^{(k)}\right) + s_{GP,res}. \tag{28}$$

### 3.2 Leveraging Second-Order Information

In computational tests with A0, the optimiser frequently entered the restoration phase, where the trust-region shrinks to improve surrogate accuracy. Although this mechanism guarantees eventual recovery, it becomes inefficient because small TRF tuning parameters ($\Delta_{min}$, $\epsilon_\theta$) may yield very small steps, slowing down the progress. Relaxing these tolerances may help, but requires manual tuning, reducing robustness.

To address this, we incorporate second-order information into the TRF algorithm, inspired by recent advances in scaling techniques [27][30] and eigenvalue-based curvature control [22][50][51]. By penalising curvature-sensitive steps using Hessian information while relaxing constraints on less sensitive directions, the method reduces reliance on parameter tuning. These Hessian-based modifications enhance robustness for poorly scaled NLPs, where standard TRF may stall despite convergence guarantees.

In the base A0 algorithm, TRSP uses $l_\infty$-norm [22]:

$$-\Delta^{(k)} \leq (x_i - x_i^{(k)}) \leq \Delta^{(k)} \; for \; i = 1, \ldots, n. \tag{29}$$

Although simple, this constraint is inefficient for problems with high anisotropy/nonlinearity or poor scaling.

A common extension introduces a scaling matrix [27][22]:

$$\left\|(x - x^{(k)})\right\|_{E^{-1}} \leq \Delta^{(k)}, \tag{30}$$

where $E = diag(u_1, \dots, u_n) \in \mathbb{R}^{n \times n}$ contains nominal scaling factors $u_i > 0$.

Inspired by the scaling approach (30) and Conn's methods of initialising trust-region size using eigenvalues for NLPs (given eigenvalues' linkage to Lagrangian multipliers [55][56]) [22], we propose integrating the local Hessian matrix $H$ into TRSP (4) (and compatibility check (6)), using a geometry-aware ellipsoidal $l_2$-norm:

$$\left\|(x - x^{(k)})\right\|_{H} \leq \Delta^{(k)}. \tag{31}$$

The Hessian is approximated during the criticality check using the *gjh* solver (Pyomo GitHub repository [57]).

**Lemma 3.1**

Consider a simple trust-region subproblem ($\min q(x) := x^T A x - 2a^T x \;\; s.t. \|x\|^2 \leq \Delta^2$), where $A \in \mathbb{R}^{n \times n}$ is a symmetric matrix, $a \in \mathbb{R}^n$ and $\Delta \in \mathbb{R}_+$. We define the corresponding Lagrangian $\mathcal{L}(x, \lambda) = x^T A x - 2a^T x - \lambda(\|x\|^2 - \Delta^2)$ with Lagrangian multiplier $\lambda$. The stationarity condition requires $\nabla_x \mathcal{L}(x^*, \lambda^*) = (A - \lambda^* I)x^* - a = 0$ at KKT point $x^*$, while the second-order sufficient optimality condition requires:

$$y^T (A - \lambda^* I) y > 0 \;\; \forall\, y \ni y^T x^* \leq 0. \tag{32}$$

Condition (32) is further satisfied with:

$$A - \lambda^* I \succcurlyeq 0. \tag{33}$$

Based on Lemma 3.1 [58][59], if $H$ is indefinite, it should be projected to a positive definite (PD) form before use in TRSP. Using an indefinite Hessian can admit multiple stationary points and compromise global convergence, particularly in non-convex settings. Projection restores standard convergence guarantees, ensuring descent directions and satisfaction of Wolfe or trust-region conditions [60][22].

### 3.2.1 Eigenvalue-based Diagonal Loading (A1 Algorithm)

Diagonal loading is a simple and robust approach to generate a PD Hessian, denoted as $H_{DL}$. For an iteration $k$, we perturb the Hessian by a scaled identity matrix so that all eigenvalues are shifted above a small positive threshold, yielding a projected PD Hessian:

$$H_{DL}^{(k)} = H^{(k)} + \tau I, \quad \tau = \max\{\varepsilon_1 - \lambda_{min}(H^{(k)}), 0\}, \tag{34}$$

where $\lambda_{min}(H^{(k)})$ is the smallest eigenvalue of the current Hessian and $\varepsilon_1 > 0$ is a small safeguard [61].

When Hessian $H^{(k)}$ is already PD, $\tau = 0$ enables the full Newton step when feasible. Otherwise, the shift raises all eigenvalues above $\varepsilon_1$, ensuring strict positivity. This eigenvalue-based adjustment prevents PD-related breakdowns in Newton or TRF methods with negligible tuning overhead. Because it only requires the smallest eigenvalue, the approach is computationally inexpensive compared to full spectral projections, making it attractive in large-scale engineering applications [62].

By employing $H_{DL}^{(k)}$, the trust-region constraint becomes:

$$\left\|(x - x^{(k)})\right\|_{H_{DL}^{(k)}} \leq \Delta^{(k)} \quad or \quad \left[x - x^{(k)}\right]^T H_{DL}^{(k)} \left[x - x^{(k)}\right] \leq \Delta^{(k)}. \tag{35}$$

We refer to this diagonal-loading modification as the A1 algorithm in the results section.

### 3.2.2 Spectral Projection Strategies (A2-A4)

While eigenvalue-based diagonal loading offers a computationally inexpensive way to reflect local curvature in trust-region formulations, incorporating the full Hessian (i.e., all eigenvalues $\lambda^{(k)}$) provides a more accurate geometric representation of the optimisation landscape. Following the use of projected Newton strategies in highly nonlinear simulations [63][50], we incorporate projected Hessians into the TRF algorithm. Two spectral-decomposition approaches are considered to regularise problematic eigenvalues [22].

**Clamped Projection (A2 Algorithm)**

In the A2 strategy, negative or zero eigenvalues are clamped to a small positive value $\varepsilon_2$ ignoring directions associated with uncertain or ill-posed curvature:

$$\lambda_C^{(k)} = \begin{cases} \varepsilon_2 \; if \; \lambda^{(k)} \leq \varepsilon_2, \\ \lambda^{(k)} \; otherwise. \end{cases} \tag{36}$$

The eigenvalues $\lambda_C^{(k)}$ define a PD clamped projected Hessian $H_C^{(k)}$ [64], used in the trust-region constraint:

$$\left[x - x^{(k)}\right]^T H_C^{(k)} \left[x - x^{(k)}\right] \leq \Delta^{(k)}. \tag{37}$$

A2 flattens the negative curvature to $\varepsilon_2$, a widely used strategy [63]. The resulting PD Hessian guarantees descent for Newton steps. Convexification ensures rigorous convergence safeguards, as standard trust-region theory [22][65] applies under mild assumptions. The drawback is loss of curvature fidelity: in strongly nonconvex regions, clamping may distort valleys or saddles, producing overly large steps toward the trust-region boundary and away from the true minimiser [51]. In practice, A2 is robust but often slow, requiring many small steps in strongly nonconvex problems [51]. Convergence can even stall when the true surface is highly curved because A2 does not exploit any negative-curvature information [50].

**Absolute Projection (A3 Algorithm)**

In the A3 strategy, negative eigenvalues are replaced by their absolute values, while zero eigenvalues are clamped to $\varepsilon_3$:

$$\lambda_A^{(k)} = \begin{cases} \varepsilon_3 \ \ if \ \left|\lambda^{(k)}\right| \le \varepsilon_3, \\ \left|\lambda^{(k)}\right| \ otherwise. \end{cases} \tag{38}$$

The resulting absolute projected Hessian $H_A^{(k)}$, defines the trust-region constraint:

$$\left[x - x^{(k)}\right]^T H_A^{(k)} \left[x - x^{(k)}\right] \le \Delta^{(k)}. \tag{39}$$

Geometrically, A3 retains curvature magnitudes (including concave directions) but not their signs, enforcing strong convexification. Algebraically, the Newton step with absolute projected Hessian is equivalent to a scaled gradient step in directions of strong negative curvature. Thus, A3 produces conservative steps, avoiding the disproportionately large moves seen in A2. A3 often accelerates convergence in ill-conditioned, noisy, or nonconvex scenarios [51]. Scaling by magnitude avoids the infinite step issue inherent to clamping. For large negative eigenvalues, A2 permits larger steps, whereas A3 produces smaller, safer steps on steep saddles. The drawback is geometric distortion: steps deviate from the original quadratic model, so the path is less geometry-aware. In trust-region terms, absolute projection acts like a vanishing-radius model (a pure gradient step) when negative curvature dominates. This slows convergence near the true minimum, and classical convergence proofs may not apply directly. In short, A3 trades model accuracy for robustness [51]: it is faster in early nonconvex phases, while A2 catches up in later convex regimes.

**Adaptive Switching (A4 Algorithm)**

Recognising that neither A2 nor A3 is uniformly better, we propose an adaptive strategy that switches between them based on surrogate accuracy. The reduction ratio $\rho$ (15) is monitored:

- If $\rho \approx 1$: the surrogate is accurate, larger steps are justified. We adopt the clamped (A2) strategy, close to Newton [51].
- If $\rho$ is small or negative: the surrogate is poor, and the trust-region shrinks. In this regime, A3 is used to enforce conservative steps [51].

The optimisation begins with the A3 strategy (appropriate for highly nonlinear and nonconvex problems) and switches adaptively as convergence proceeds. The complete switching logic is illustrated in Figure 3. During optimisation using the A4 algorithm (Figure 3), if surrogate agreement is strong (i.e., $f - type$ or $\rho^{(k)} \ge \eta_2$ in $\theta - type$/restoration step), A2 is applied. Otherwise (i.e., rejected step or $\rho^{(k)} < \eta_2$ in $\theta - type$/restoration step), the method switches back to A3 for robustness [50][51].

This adaptive rule exploits curvature when reliable (A2) and ensures safety when not (A3).

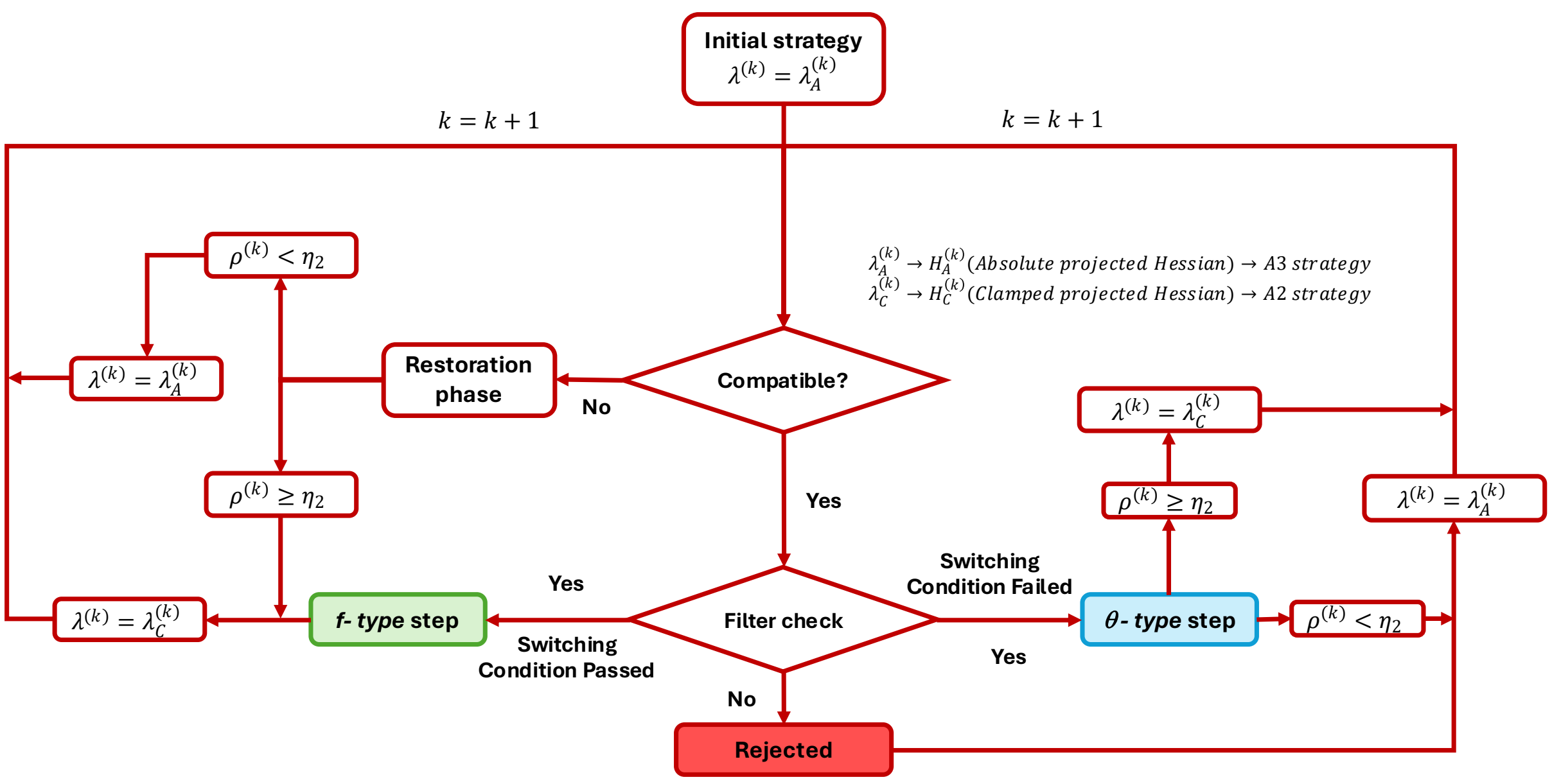

***Figure 3.*** *A4 algorithm: adaptive projected-Hessian filtering embedded within the TRF framework. The scheme shows switching between clamped (A2) and absolute (A3) Hessian projections based on the surrogate accuracy and trust-region reduction ratio. When surrogate agreement is high, A2 is selected for curvature exploitation; when poor, A3 enforces conservative steps. The adaptive strategy balances robustness and efficiency across nonlinear grey-box problems.*

### 3.3 Dimensionality Reduction Strategy (A5) and Algorithms Summary

A dimensionality reduction strategy has previously been implemented in the simplified TRF method (A5), an extension of the classical TRF method (A0), with a TS surrogate using a global basis [27]. We use the A5 method as another benchmark in addition to the A0 algorithm. In A5, variables $x^T = [w^T, y^T, z^T]$ (as in A0) are partitioned into decision variables $u$ and remaining variables $d$:

$$x^T = [u^T, d^T].$$

The trust-region constraint (30) is applied only to $u$, reducing dimensionality and simplifying scaling. On the other hand, A5 [27] uses explicit derivative information and prior classification of decision variables. In A1-A4, near-zero eigenvalues naturally reduce dimensionality in the TRSP by relaxing steps in corresponding directions. Compared to A0 and A5, Hessian-based A1-A4 methods provide stability in poorly scaled or anisotropic spaces. It should be noted that eigenvalue modifications in A1-A4 do not strictly satisfy the classical convergence assumptions of model-based trust-region methods. Thus, these methods should be viewed as practical heuristics: relaxing strict theory in exchange for robustness in challenging grey-box optimisation.

Table 1 summarises the six TRF variants (A0-A5) used/introduced in this work. A0 serves as the baseline classical TRF benchmark, while A1-A4 introduce Hessian regularisation and projection techniques, ranging from simple diagonal loading to adaptive spectral projections, to improve robustness and reduce sensitivity to tuning parameters. The simplified TRF (A5) [10][27] applies a reduced-dimensional trust-region constraint to manually selected decision variables. Together, these variants illustrate the systematic evolution from the classical TRF methods to the proposed Hessian-informed variants, achieving greater numerical stability and efficiency.

***Table 1.*** *Summary and comparison of the six TRF variants (A0-A5) introduced/used in this study. Each variant differs in its trust-region (TR) constraint formulation and Hessian regularisation. A0 represents classical TRF, A1-A4 introduce Hessian-based TRF strategies for improved numerical stability, and A5 corresponds to the simplified TRF method that applies a reduced-dimensional TR constraint. The "Distinct Behaviour" column highlights key algorithmic features and performance characteristics.*

| Variant | TR Constraint Formulation | Hessian Regularisation | Hessian Projection Type | Distinct Behaviour |
|---|---|---|---|---|
| A0 | $\left\|(x - x^{(k)})\right\|_{\infty} \leq \Delta^{(k)}$ | None | None | Classical TRF; sensitive to parameter tuning and poorly scaled problems. |
| A1 | $\left\|(x - x^{(k)})\right\|_{H_{DL}^{(k)}} \leq \Delta^{(k)}$ | Diagonal Loading | Diagonal Loading | Perturbs the Hessian by a scaled identity matrix; inexpensive and stable for large systems. |
| A2 | $\left\|(x - x^{(k)})\right\|_{H_{C}^{(k)}} \leq \Delta^{(k)}$ | Spectral | Clamped Eigenvalues | Clamps negative/zero eigenvalues to small positive values; guarantees descent but may distort curvature in nonconvex regions. |
| A3 | $\left\|(x - x^{(k)})\right\|_{H_{A}^{(k)}} \leq \Delta^{(k)}$ | Spectral | Absolute Eigenvalues | Uses eigenvalue magnitudes to preserve curvature scale; robust for ill-conditioned or nonconvex problems. |
| A4 | $\left\|(x - x^{(k)})\right\|_{H_{A/C}^{(k)}} \leq \Delta^{(k)}$ | Spectral | Clamped/Absolute Eigenvalues | Adaptively switches between A2 and A3 based on surrogate accuracy; balances curvature use and step conservatism. |
| A5 | $\left\|(u - u^{(k)})\right\|_{E^{-1}} \leq \Delta^{(k)}$ | None | None | Simplified TRF; applies trust-region constraint to manually selected decision variables and requires explicit derivatives at the current point. |

### 3.4 Termination Criteria

In addition to the original termination conditions in A0 (which involve $\epsilon_{\theta}$, $\epsilon_{\chi}$, $\epsilon_{\Delta}$, and $\Delta_{min}$), we introduce a residual-based optimality condition that enhances robustness across algorithmic variants and diverse problem sets. Since the Hessian is regularised to be PD, all usual assumptions of trust-region convergence theory remain intact [1]. In particular, if the convexified TRSP admits $r^{(k)} = 0$ as its minimiser, then the gradient of the TRSP vanishes in the subspace of active constraints. By first-order consistency, this implies that the KKT conditions hold both for the subproblem and, up to linearisation error, for the true problem.

Accordingly, we require termination whenever both feasibility and residual optimality are satisfied:

$$\theta^{(k)} \leq \epsilon_{\theta} \ \ and \ \ \left\|r^{(k)}\right\| \leq \epsilon_{r}.$$

This dual criterion guarantees that the iterate is simultaneously nearly feasible and sufficiently close to a stationary point [65][9]. The automatic scaling properties in the TRF variants (A1-A4) further reinforce this approach, since they align the residual measure with the local geometry, making $\epsilon_{r}$ easier to tune.

Thus, residual-based termination is a practical addition rather than a replacement of the classical criteria, preserving convergence guarantees while ensuring more consistent practical performance across nonlinear grey-box optimisation problems.

# 4. Implementation

The TRF methods (A0-A4) are implemented in a Python-based Pyomo modelling environment [57], which provides a flexible open-source platform for developing optimisation algorithms. All subproblems are solved using IPOPT 3.14.13 with linear solver ma27 [7]. Alternative NLP solvers can be used as needed. Derivatives and Hessians required for the criticality check and the modified trust-region constraints are obtained from the AMPL pseudo-solver *gjh*. Numerical experiments are executed on a Microsoft Windows 11 system with an Intel(R) Core(TM) i7-1165G7 processor (8 cores) and 16 GB RAM. All runs are executed in serial mode (no explicit parallelisation) using the default single-threaded settings of the solver. All codes for A0-A4 are available in our GitHub repository [66], as referenced in the Data Availability Statement.

# 5. Numerical Results

Since no benchmark library exists for grey-box optimisation, we compiled a test set of 25 NLPs with black-box components (simulation experiments [67] and published literature [68]). In addition, five engineering optimisation problems were formulated. In these problems, one or more (up to four) nonlinear constraints or parts of the objective function are treated as black-box, while retaining the remaining part of the model as glass-box. All Pyomo codes for each problem are available in the GitHub repository [66].

The performance of the proposed trust-region filter (TRF) algorithms was evaluated through a three-tier analysis designed to assess robustness, efficiency, and scalability. First (section 5.1), the new Hessian-informed TRF variants (A1–A4) were benchmarked against the baseline classical TRF (A0) and the simplified dimensionality-reduction variant (A5) to quantify the impact of curvature information and surrogate selection. Second (section 5.2), the best-performing grey-box Hessian-based TRF algorithms were compared to the state-of-the-art DFO solvers (COBYLA, Py-BOBYQA, NEWUOA, Nelder-Mead, PyNomad and COBYQA), validating their competitiveness and the advantages of exploiting glass-box model structure. Finally (section 5.3), the proposed algorithms were applied to engineering case studies, demonstrating their practical performance and scalability in grey-box environments.

## 5.1 Benchmarking Against Available Grey-box TRF Algorithms

In this section, the proposed Hessian-informed TRF variants (A1–A4) are first compared with the base TRF algorithm (A0) across six surrogate types (linear, standard quadratic, simplified quadratic, GP, TS and hybrid). Following that, for completeness, we also compare the proposed methods with the simplified TRF algorithm (A5) under identical settings.

**Comparison with A0**

Unlike A0 (which requires tuning of several critical parameters: $\epsilon_\theta$, $\epsilon_\chi$, $\epsilon_\Delta$, $\Delta^{(0)}$ and $\Delta_{min}$), the Hessian-based variants rely on minimal tuning ($\epsilon_\theta$ and $\epsilon_r$) or none at all, owing to automatic step scaling by the projected Hessian.

Performance is summarised in Table 2 as the percentage of problems solved, revealing distinct surrogate performance clusters. Higher-fidelity surrogates (GP and TS) demonstrated superior robustness, solving 92-100% of the test set. All Hessian-based variants (A1-A4) paired with GP or TS achieved 100% success. Polynomial surrogates (linear, quadratic, simplified quadratic) performed moderately, solving 72-84% of the problems, with the best results obtained using A3 (84%), followed by A2 and A0 (80-84%), A4 (72-84%), and A1 (72-80%). The hybrid surrogate was least effective, solving fewer than 16% of problems, suggesting poor compatibility with the TRF framework. Across all surrogate choices, the spectrally projected Hessian variants A2 and A3 solved the maximum number of problems.

***Table 2.*** *Success rate of algorithmic variants (A0–A4) combined with different surrogate models across 25 benchmark problems. Rows correspond to TRF variants and columns to surrogate types (linear, quadratic, simplified quadratic, Gaussian processes - GP, Taylor series - TS, hybrid). Values denote the fraction of problems solved successfully. High-fidelity surrogates (GP and TS) achieved 92-100 % success compared with 72-84 % for low-fidelity surrogates.*

| | **Linear (L)** | **Standard Quadratic (Q)** | **Simplified Quadratic (SQ)** | **Gaussian Process (GP)** | **Taylor Series (TS)** | **Hybrid (H)** | **Variants Average Performance** |
|---|---|---|---|---|---|---|---|
| **A0** | 0.80 | 0.84 | 0.80 | 0.92 | 0.92 | 0.12 | 0.733 |
| **A1** | 0.76 | 0.80 | 0.72 | 1 | 1 | 0.12 | 0.733 |
| **A2** | 0.84 | 0.84 | 0.80 | 1 | 1 | 0.16 | 0.773 |
| **A3** | 0.84 | 0.84 | 0.84 | 1 | 1 | 0.16 | **0.780** |
| **A4** | 0.84 | 0.80 | 0.72 | 1 | 1 | 0.12 | 0.747 |
| **Surrogates Average Performance** | 0.82 | 0.82 | 0.78 | **0.98** | **0.98** | 0.14 | |

A comparative analysis of algorithmic variants and surrogate models was conducted using three performance metrics: number of iterations, computational time (s), and external function evaluations. Computational times represent wall-clock runtime measured in Python. Figure 4, Figure 5 and Figure 6 consist of five subplots, each corresponding to one of the TRF variants (A0-A4). The x-axis represents iterations (Figure 4), external evaluations (Figure 5) and computational time (Figure 6) on a logarithmic scale, while the y-axis indicates the percentage of problems solved. Each curve corresponds to a different surrogate. A steeper rise in the curve indicates more problems solved with a lower budget, while a higher final plateau reflects greater robustness. Flat regions signal stagnation, where additional iterations or evaluations failed to yield further progress.

Figure 4 highlights the number of iterations required to solve the benchmarks across algorithmic-surrogate configurations. Hessian projection strategies (A1-A4) consistently improved upon the base algorithm A0 by stabilising the TRSP and exploiting curvature information. Without these mechanisms, A0 often entered the restoration phase or had successive steps rejected, gradually shrinking the trust-region radius. Progress could only be forced by resetting parameters such as $\Delta_{min}$ and $\epsilon_{\theta}$ to larger values, yet the iteration count remained high. In contrast, projected Hessians (A1-A4) automatically scaled and penalised steps, reducing sensitivity to tuning parameters and rarely requiring restoration. Diagonal loading (A1) modestly accelerated convergence, while eigenvalue-based projections (A2, A3) offered significant improvements: for example, under A2 (clamped Hessian) with TS and GP surrogates, nearly all problems were solved in fewer than six iterations, compared to tens of iterations for A0 with the same surrogate. A3 (absolute Hessian) often matched or outperformed A2, consistent with the theory that absolute filtering avoids pitfalls of clamping. The adaptive switching strategy (A4) combined the strengths of both, quickly leveraging whichever projection yielded a more reliable step. Overall, all Hessian-modified variants (A1-A4) achieved faster convergence than A0, with larger, more reliable steps enabled by enforcing a PD subproblem Hessian.

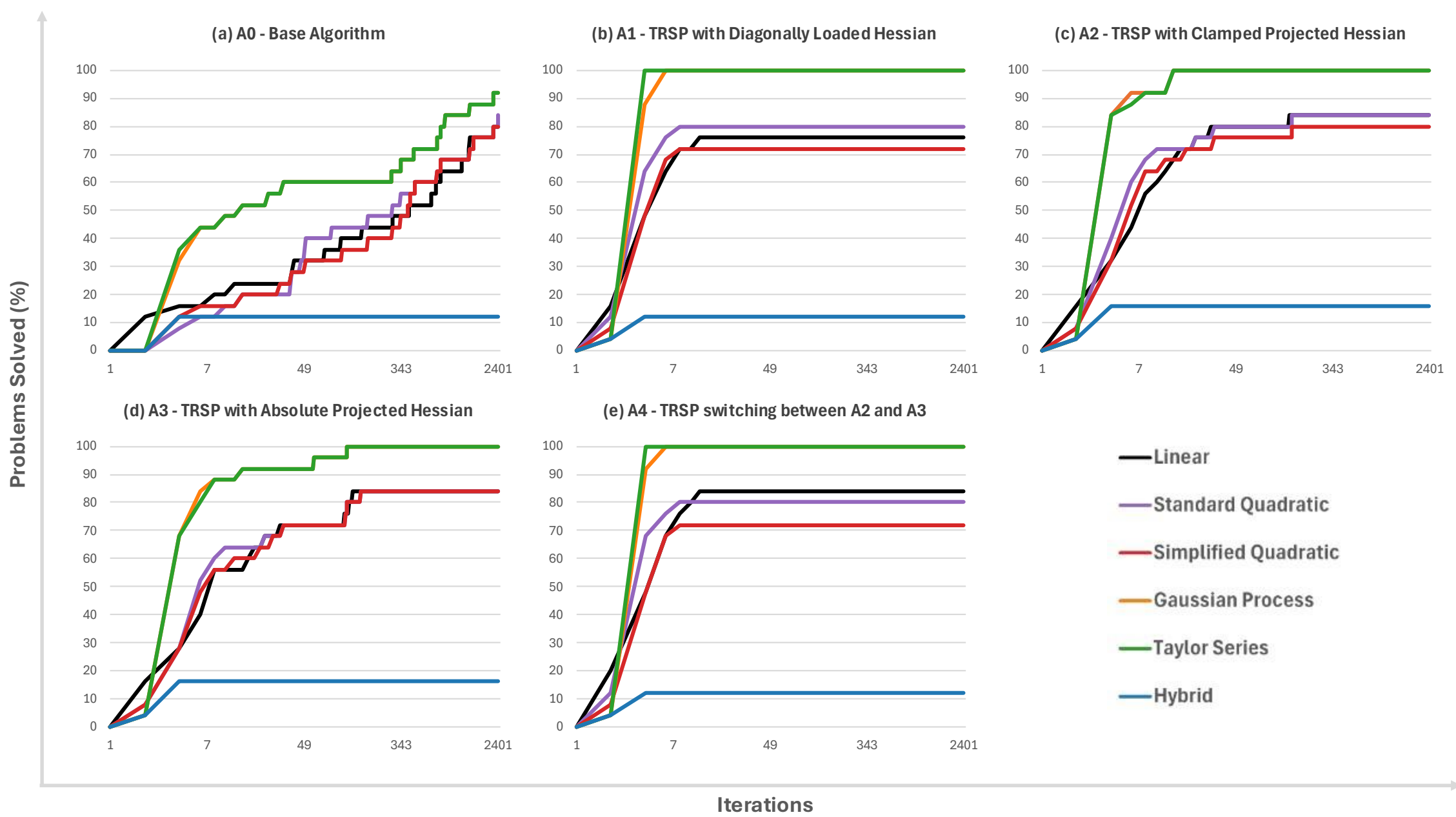

***Figure 4.*** *The number of iterations versus the percentage of problems solved successfully. Hessian-informed variants (A1-A4) achieve faster convergence than the baseline A0 across all surrogate types, particularly when using high-fidelity surrogates (Gaussian processes and Taylor series).*

Iteration counts alone do not capture computational cost. Figure 5 and Figure 6 show external evaluations and CPU time, respectively, which reflect surrogate training overhead. GP surrogates, though highly accurate, require more time to train than TS or polynomial surrogates, while quadratic surrogates generally require more evaluations. Nonetheless, the trends mirror those seen in iteration performance: TS and GP surrogates paired with Hessian-projected methods achieved rapid and reliable convergence, while polynomial surrogates improved noticeably when paired with A1-A4 compared to A0.

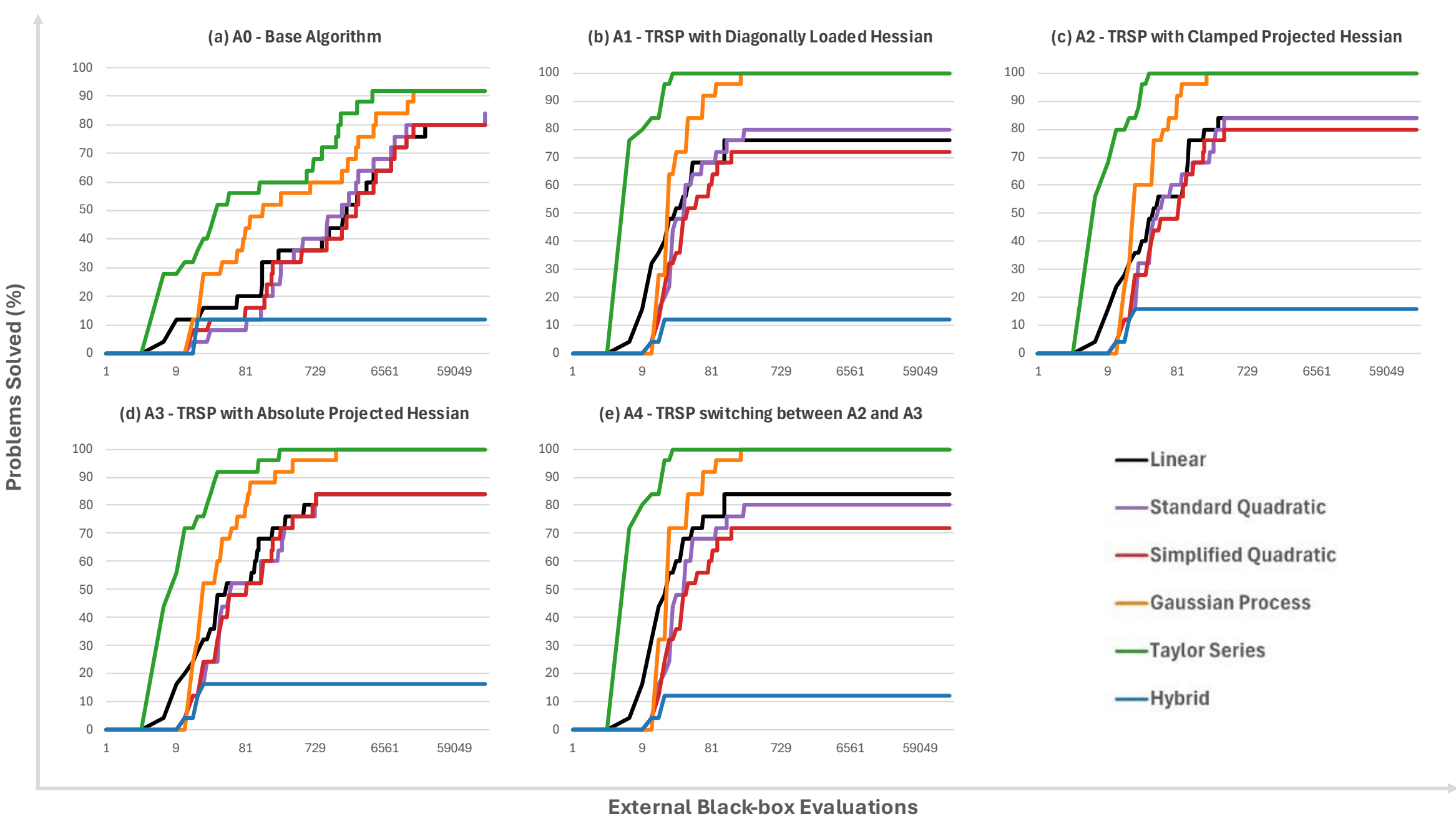

***Figure 5.*** *External black-box function evaluations versus percentage of problems solved successfully. Projected-Hessian variants (A1-A4) reduce the number of expensive black-box calls by up to an order of magnitude relative to the classical TRF (A0), demonstrating improved computational efficiency.*

TS surrogates proved most effective overall, solving nearly all problems rapidly with any enhanced TRF variant (A1-A4). GP surrogates also paired strongly with Hessian projection, particularly A1, A2, and A4, solving 100% of problems within 204 evaluations or 10-17 seconds. Polynomial surrogates (linear, quadratic, simplified quadratic) solved 72-84% of problems, requiring more steps (124-354 evaluations). However, they still benefited substantially from Hessian projection: A1-A4 with polynomials consistently outperformed the A0-polynomial combinations. The hybrid surrogate consistently underperformed, indicating that naive model combinations do not produce reliable reduced models in the TRF setting.

While GP surrogates required higher training time, they offered important advantages; they eliminated the need for derivative information, captured black-box structure effectively, and halved the number of external black-box calls relative to polynomial surrogates [29][30]. As such, GP-based TRF methods remain highly attractive when derivatives are unavailable or prohibitively expensive. Moreover, GP surrogates have well-established convergence properties in derivative-free optimisation [10], making them a robust choice.

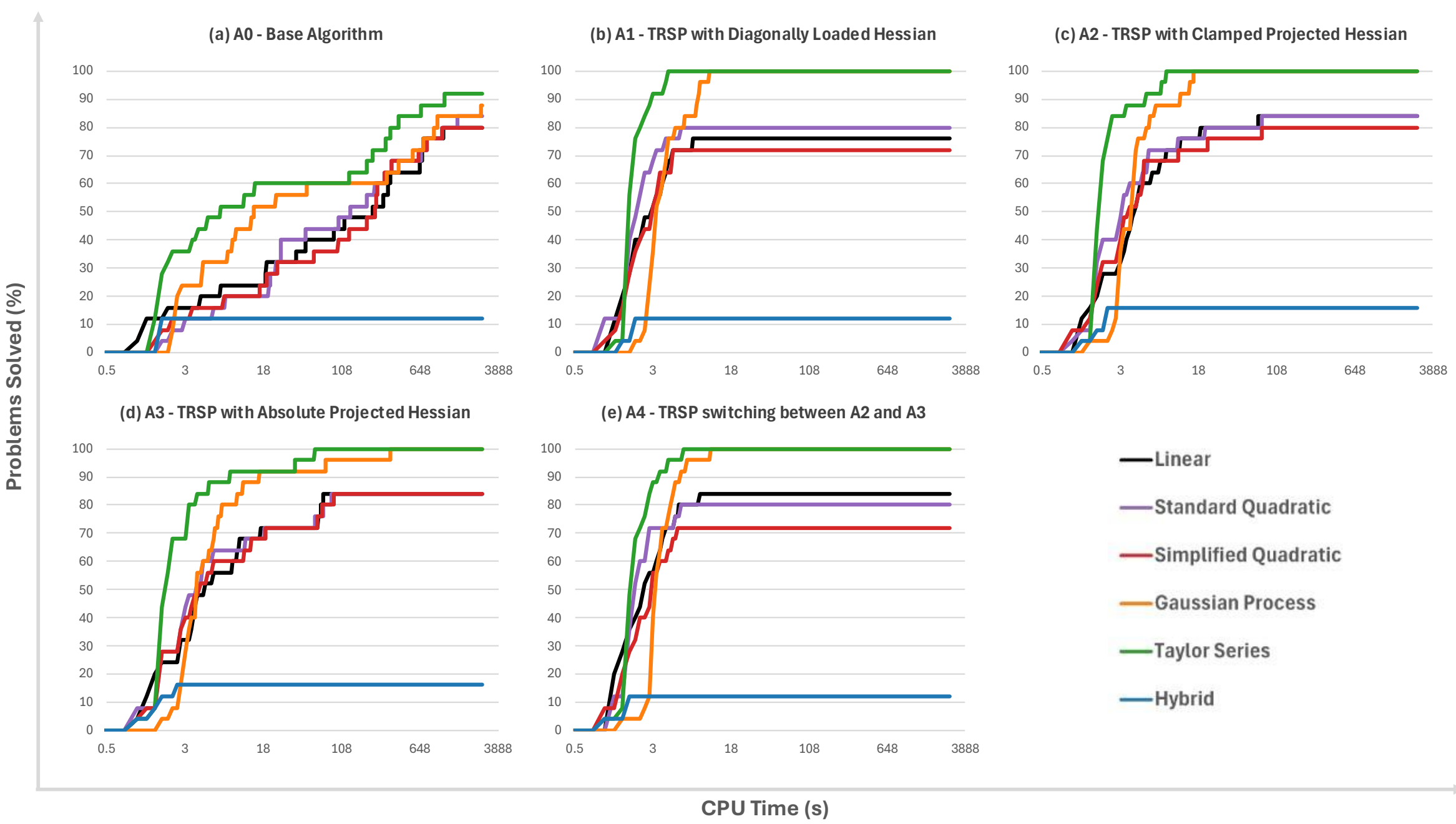

***Figure 6.*** *Computational time versus percentage of problems solved successfully for different TRF variants and surrogate models. Hessian-regularised algorithms achieve lower total runtime, with the adaptive variant A4 offering the best trade-off between accuracy and speed.*

To directly compare the algorithmic variants independent of surrogate model effects, we constructed Dolan-Moré performance profiles [69] using CPU time (s) as a performance metric, representative of both external black-box evaluation cost and iteration trends for the problem set. For each algorithm $s$ and problem $p$, the performance ratio is defined as

$$r_{p,s} = \frac{t_{p,s}}{min\{t_{p,s}: s \in S\}}, \tag{40}$$

where $t_{p,s}$ is the CPU time required by algorithm $s$ to solve problem $p$, and $S$ is the set of all algorithms tested. The performance profile for algorithm $s$ is then expressed as

$$\rho_s(\alpha) = \frac{1}{|P|} size\{p \in P: r_{p,s} \leq \alpha\}, \tag{41}$$

where $|P|$ denotes the cardinality of problem set $P$. Thus, $\rho_s(\alpha)$ represents the fraction of problems solved by algorithm $s$ within a factor $\alpha$ of the best solver. Figure 7 compares the five algorithmic variants (A0-A4) for each surrogate type. The Hessian-based variants (A1-A4) outperform the base algorithm A0 with all surrogates, achieving faster convergence across the benchmark set. Among Hessian variants, A1 and A4 show the most uniform and improved performance across the benchmark set.

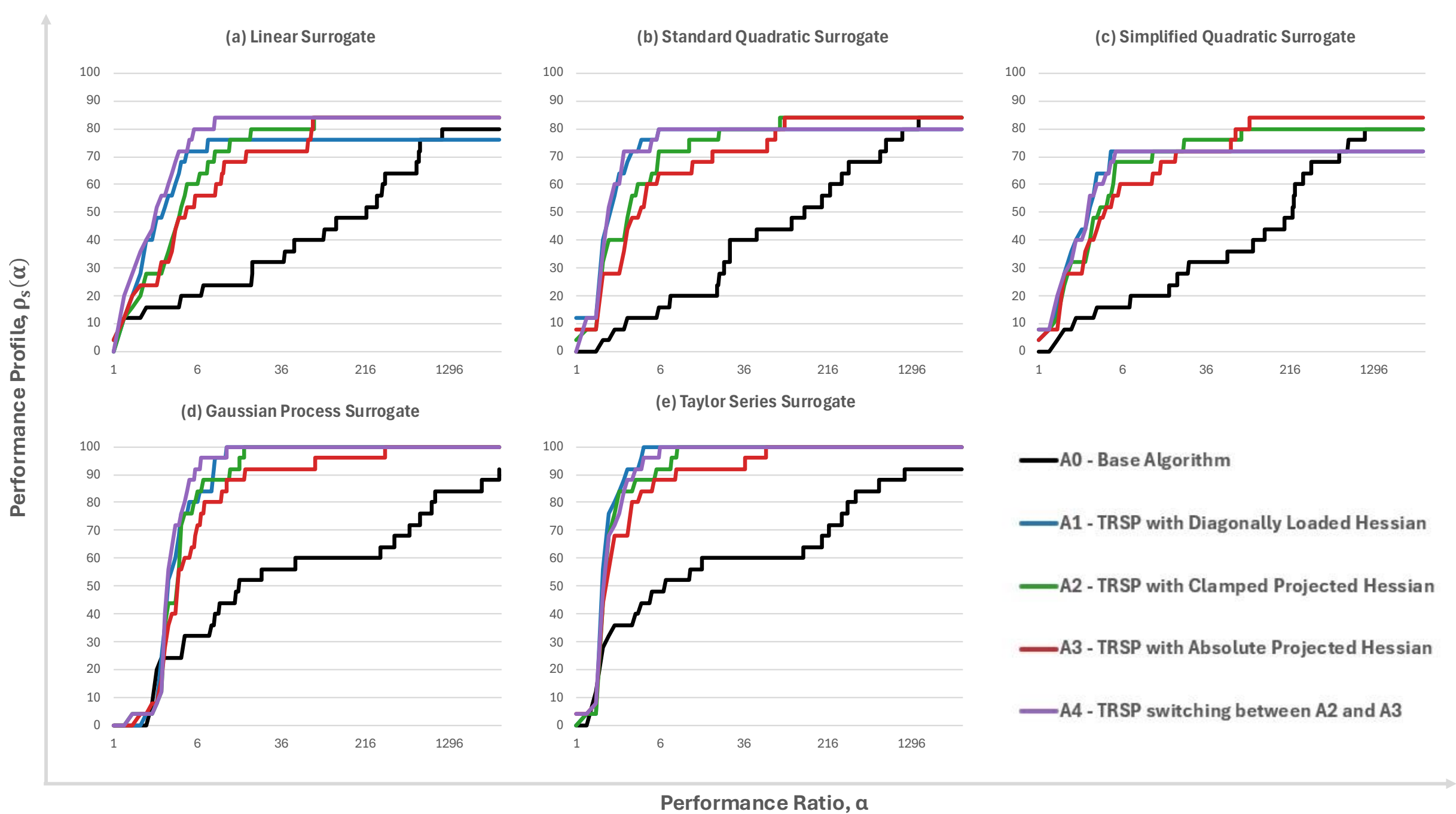

***Figure 7.*** *Dolan-Moré performance profiles comparing the five TRF algorithmic variants (A0-A4) across all surrogate models using CPU time as the performance metric. The x-axis shows the performance ratio $\alpha$ (logarithmic scale), and the y-axis indicates the performance profile (percentage/fraction of problems solved) within a factor $\boldsymbol{\alpha}$ of the best solver. Higher curves correspond to faster and more robust algorithms. The Hessian-based variants (A1-A4) consistently outperform the base TRF (A0), with A1 and A4 demonstrating the best performance across most surrogate types.*

## Comparison with A5

For completeness, the simplified TRF method (A5) was tested under identical initial settings (Figure 8). As described in Section 3.2.2, A5 applies the trust-region constraint only to a subset of manually pre-selected decision variables $u$, using explicit derivatives. A5 solved 80% of the problems within 500 iterations, converging faster than A0 (92% solved in 2194 iterations) but falling well short of Hessian-based methods A1-A4 paired with TS, which solved 100% of problems within at most 116 iterations and 234 external evaluations. Importantly, A5 is tailored for process flowsheet optimisation, whereas A0-A4 are designed as generic solvers; this distinction is reflected in their relative performance.

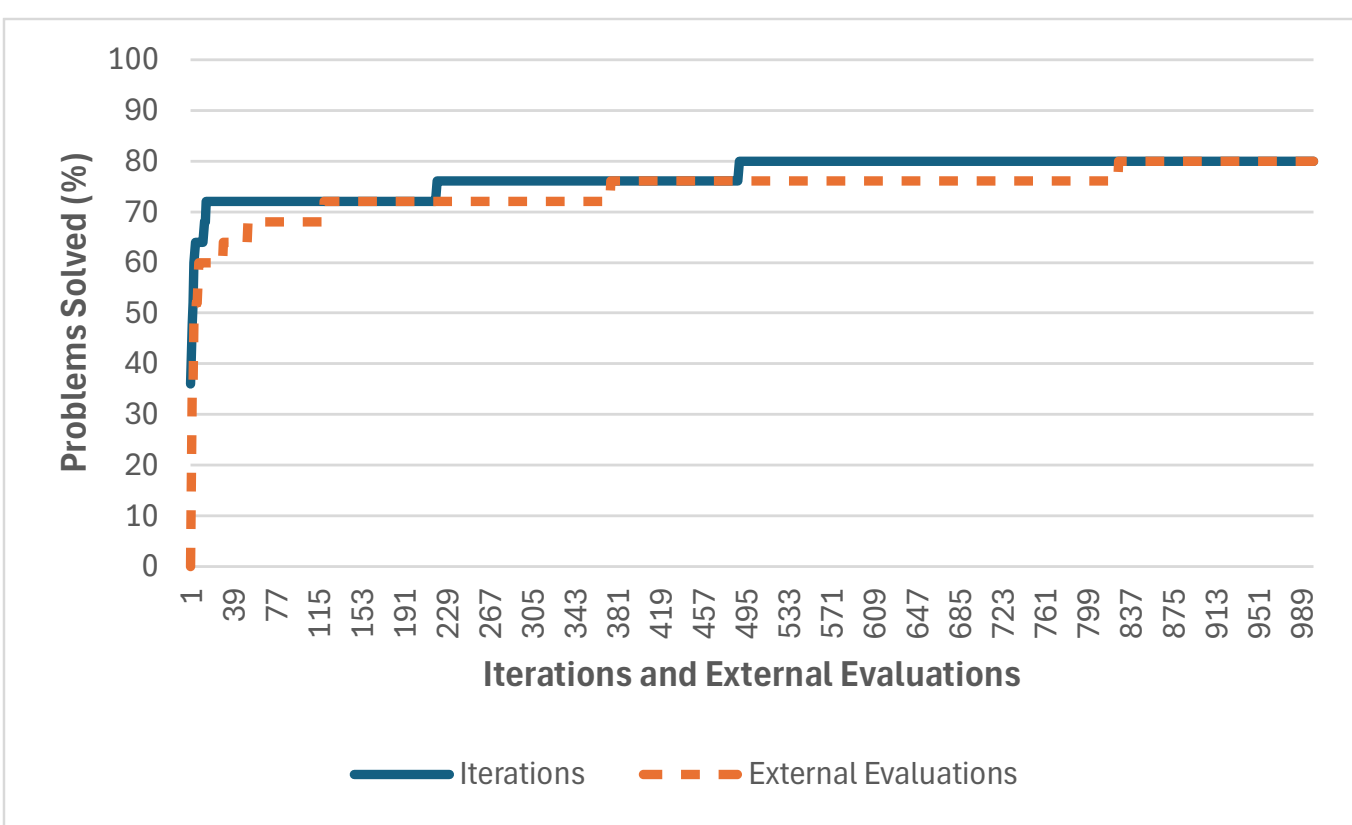

***Figure 8.*** *Performance of the simplified TRF algorithm (A5) on the grey-box test suite. While A5 provides a reduced-dimensional formulation, it underperforms compared to Hessian-based variants when derivatives are unavailable, confirming the benefit of curvature-aware scaling.*

## 5.2 Benchmarking Against State-of-the-Art Black-box DFO Solvers

Among derivative-free surrogate and variant configurations, GP paired with A1 and A4 showed the best overall performance. To further validate their robustness and generality, the two grey-box algorithms were benchmarked against several state-of-the-art DFO solvers: trust-region interpolation-based (COBYLA, Py-BOBYQA, NEWUOA, COBYQA), directional direct-search (PyNomad/NOMAD: a generalised pattern search and mesh adaptive direct search) and simplicial direct-search methods (Nelder–Mead), representing the most widely used black-box solvers in the literature and open-source implementations.

All problems in the 25-case test set were reformulated as pure black-box models, where both objective and constraint functions were evaluated externally. This setup provided a fair comparison between the proposed grey-box TRF algorithms, which leverage glass-box model derivatives and curvature information, and traditional black-box DFO solvers, which rely solely on sampled evaluations. For unconstrained or box-constrained methods (e.g., NEWUOA, Nelder–Mead, and Py-BOBYQA), constraints were incorporated via penalty formulations, while COBYLA, COBYQA and PyNomad were applied in their native constrained forms.

Performance was assessed in terms of objective accuracy and black-box function evaluations, and summarised through Dolan–Moré performance profiles (Figure 9) constructed using equations (40) and (41). The profiles (based on black-box evaluations) demonstrate that A1-GP and A4-GP consistently outperform all classical DFO solvers across the benchmark set, achieving superior performance. Furthermore, both grey-box TRF variants achieved comparable or better optima when compared to DFO results. These results highlight the fundamental advantage of incorporating curvature information and partial model structure within the TRF framework, yielding faster convergence and stronger scalability than purely data-driven DFO approaches. Such performance is particularly valuable in process optimisation settings, where some process units can be modelled analytically while others rely on costly black-box simulations.

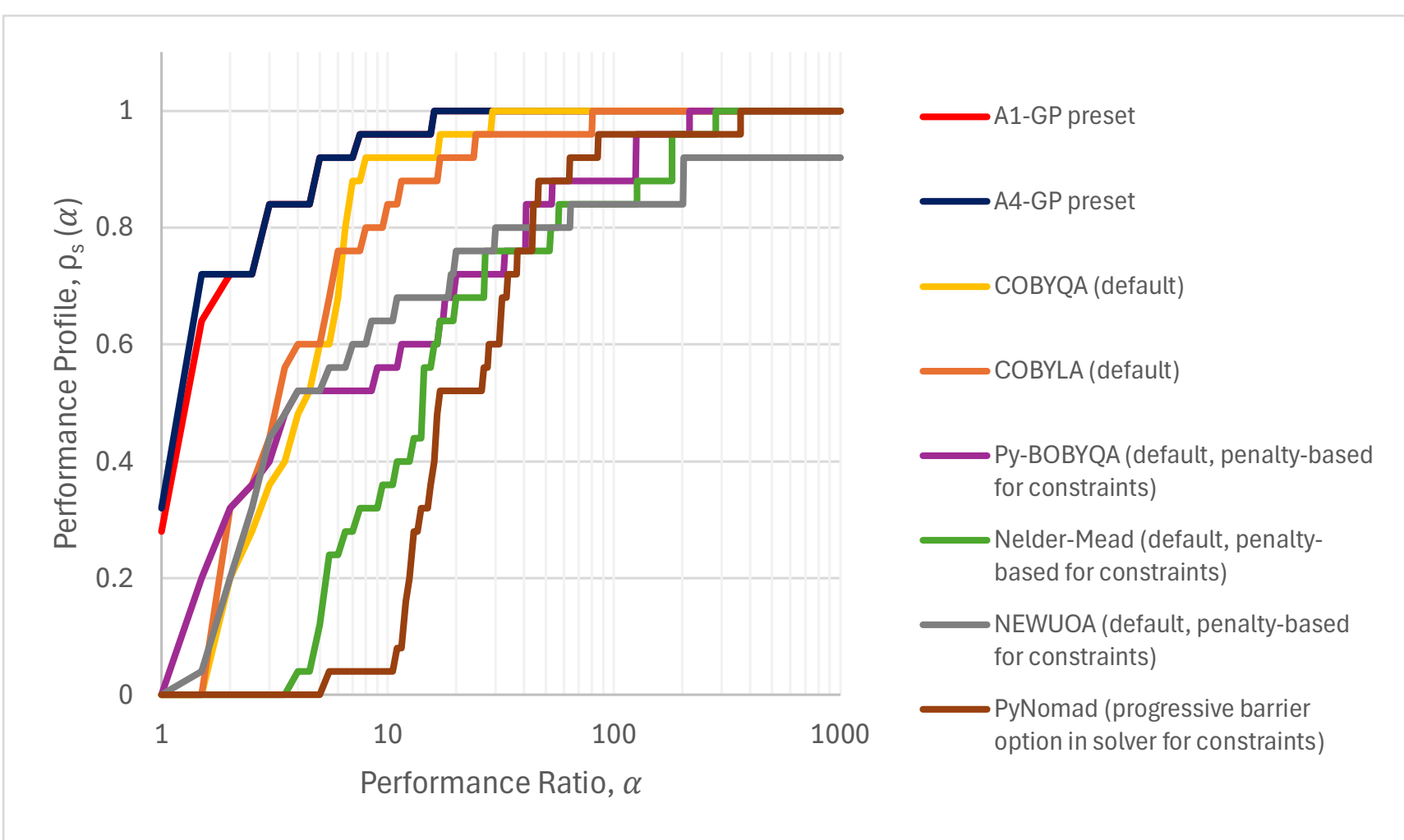

***Figure 9.*** *Performance profiles (Dolan–Moré) comparing the proposed grey-box TRF algorithms (A1–GP and A4–GP) with leading black-box DFO solvers such as COBYLA, NEWUOA, COBYQA, Py-BOBYQA, PyNomad, and Nelder–Mead, across the 25-problem benchmark set. The x-axis shows the performance ratio α (based on black-box evaluation count) on a logarithmic scale, and the y-axis represents the performance profile (percentage/fraction of problems solved) within a factor α of the best solver. The proposed A1-GP and A4-GP methods achieve the best performance, demonstrating the benefits of curvature-guided grey-box optimisation.*

## 5.3 Engineering Case Studies

In addition to the benchmark problem set, the TRF algorithms (A0-A4) were evaluated on process systems (liquid-liquid extraction and an alkylation process) and representative engineering problems (Himmelblau's problem, pressure vessel design, and spring design). These case studies were chosen to reflect a spectrum of real-world complexity, including nonlinear objectives, mixed mechanistic-empirical (grey-box) models, constrained operation, multimodality, and ill-conditioned behaviour.

### Himmelblau's problem

This well-known mechanical engineering optimisation problem $(n_w, n_y, n_z) = (3,2,5)$ [70] shares several features with industrial problems: nonlinear objectives, nonlinear equality and inequality constraints, variable bounds, multiple local minima, and strong interactions among decision variables. Such characteristics are frequently encountered in systems governed by thermodynamics, reaction kinetics, and process or mass balances. The Hessian matrix in this problem often becomes poorly conditioned, particularly near saddle points, which makes pure Newton methods unreliable without regularisation or trust-region safeguards. A1-A4 methods maintain PD of the subproblem Hessian, ensuring robust convergence to valid local optima.

Our results (Figure 10) show that the base algorithm A0 required the most iterations to converge, stalling near stationary local minima. Projection-based methods significantly improved performance: A2 (clamped Hessian) and A3 (absolute Hessian) converged more efficiently. The best performance was consistently achieved by A4 (adaptive switching) and A1 (diagonal loading), which outperformed all other variants across iterations, computational time, and external evaluations.

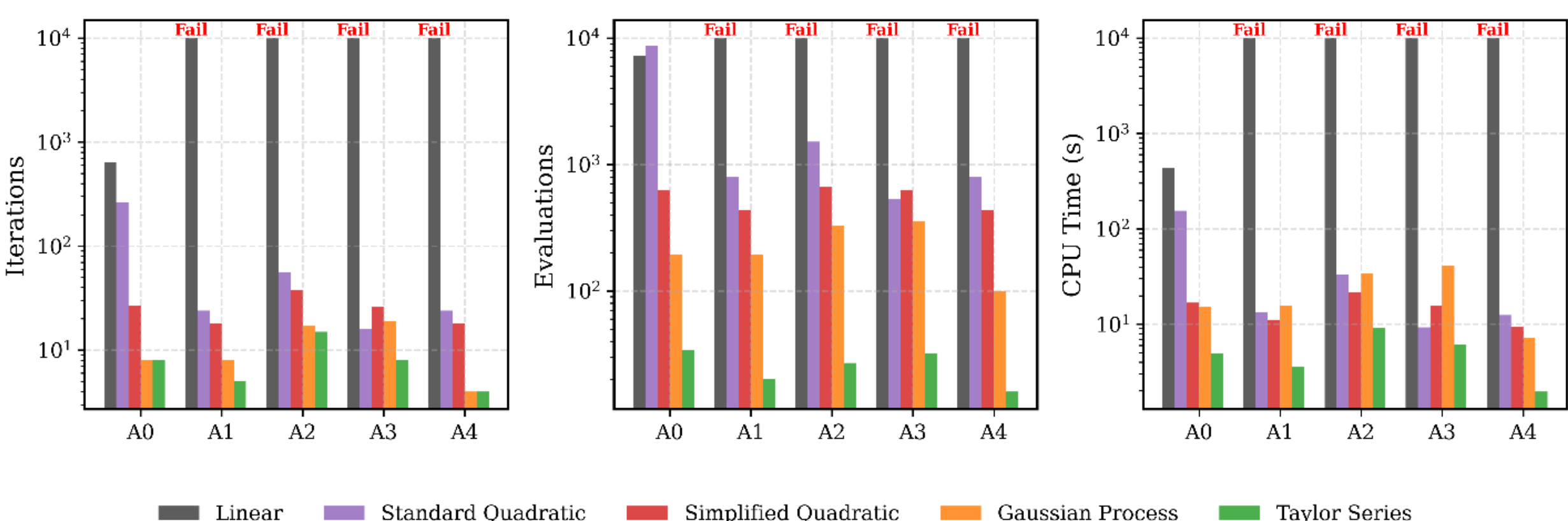

***Figure 10.*** *Comparative performance of TRF variants (A0-A4) on the Himmelblau grey-box optimisation problem. Each bar corresponds to a surrogate model. The bar with "Fail" indicated at the top shows that the specific surrogate (linear when paired with A1-A4) couldn't converge. The panels show the number of iterations (left), external black-box evaluations (middle), and computational time (right) on the y-axis (logarithmic scale). All variants converge to the same optimum; however, A1 and A4 require fewer iterations and black-box evaluations compared to other variants, indicating faster progress toward feasibility and optimality.*

## Liquid-liquid extraction column

The grey-box optimisation model of a liquid-liquid extraction column $(n_w, n_y, n_z) = (2,1,2)$ couples a mechanistic plug flow model with linear equilibrium behaviour and an empirical correlation for the number of transfer units [71]. The coupling introduces stiffness, but the Hessian remains moderately well-behaved. Nevertheless, overly aggressive trust-region steps can destabilise the solution. Within our TRF framework, Hessian information ensured well-conditioned curvature throughout the iterations. Algorithms A1-A4 consistently outperformed the base A0 method, converging with fewer iterations, requiring fewer external black-box evaluations, and achieving reduced computational time (Figure 11). These results highlight the benefit of embedding second-order curvature information into the TRF method when optimising hybrid mechanistic-empirical process models.

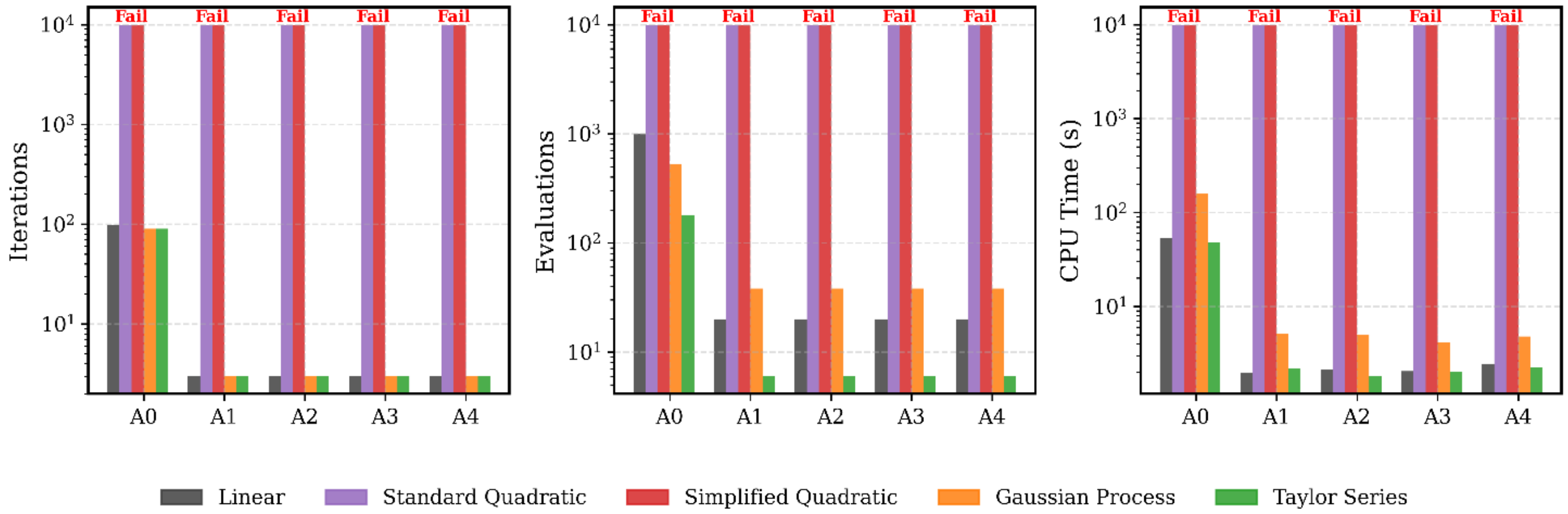

***Figure 11.*** *Comparative performance of TRF variants (A0-A4) on the liquid-liquid extraction column optimisation problem. Each bar corresponds to a surrogate model. The bars with "Fail" indicated at the top show that the specific surrogates (standard and simplified quadratic) couldn't converge. The panels show the number of iterations (left), external black-box evaluations (middle), and computational time (right) on the y-axis (logarithmic scale). The Hessian-based variants (A1-A4) outperform the classical TRF (A0) method, highlighting their effectiveness for nonlinear grey-box process optimisation.*

## Pressure vessel design

The pressure vessel design optimisation problem $(n_w, n_y, n_z) = (4,1,0)$ aims to minimise material cost subject to multiple nonlinear stress and dimensional constraints [72]. The Hessian in this problem is inherently ill-conditioned, with some directions nearly flat, destabilising gradient-based methods. By incorporating Hessian projections within the TRF framework, we regularised the curvature, preventing numerical breakdown and enabling stable progress. A1-A4 consistently solved the problem in significantly fewer iterations and with fewer external evaluations compared to the base algorithm A0 (Figure 12), highlighting the advantages of second-order information in handling ill-conditioned engineering design problems.

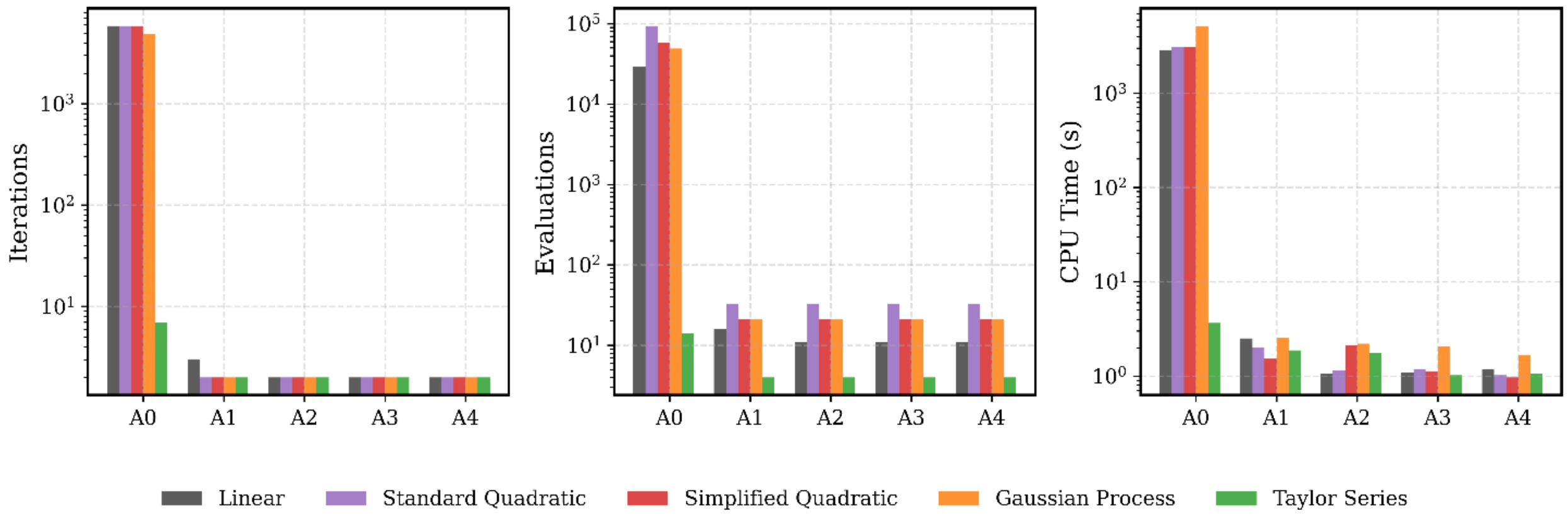

***Figure 12.*** *Comparative performance of TRF variants (A0-A4) on the pressure vessel design optimisation problem. Each bar corresponds to a surrogate model. The panels show the number of iterations (left), external black-box evaluations (middle), and computational time (right) on the y-axis (logarithmic scale). Curvature-informed variants (A1-A4) reduced iterations and external black-box evaluations by roughly an order of magnitude relative to the baseline (A0) algorithm.*

## Alkylation process

The alkylation process optimisation problem $(n_w, n_y, n_z) = (1,1,9)$ combines mechanistic mass balances with regression-based models for yield and quality variables [73]. The Hessian in this problem exhibits moderate variability with positive, negative and zero eigenvalues, making step reliability sensitive to the choice of local approximation. Within our TRF framework, trust-region steps guided by curvature information and the filter acceptance criterion helped avoid unreliable linearisations. As a result, A1-A4 consistently converged in fewer iterations and with lower computational burden compared to the A0 algorithm (Figure 13).

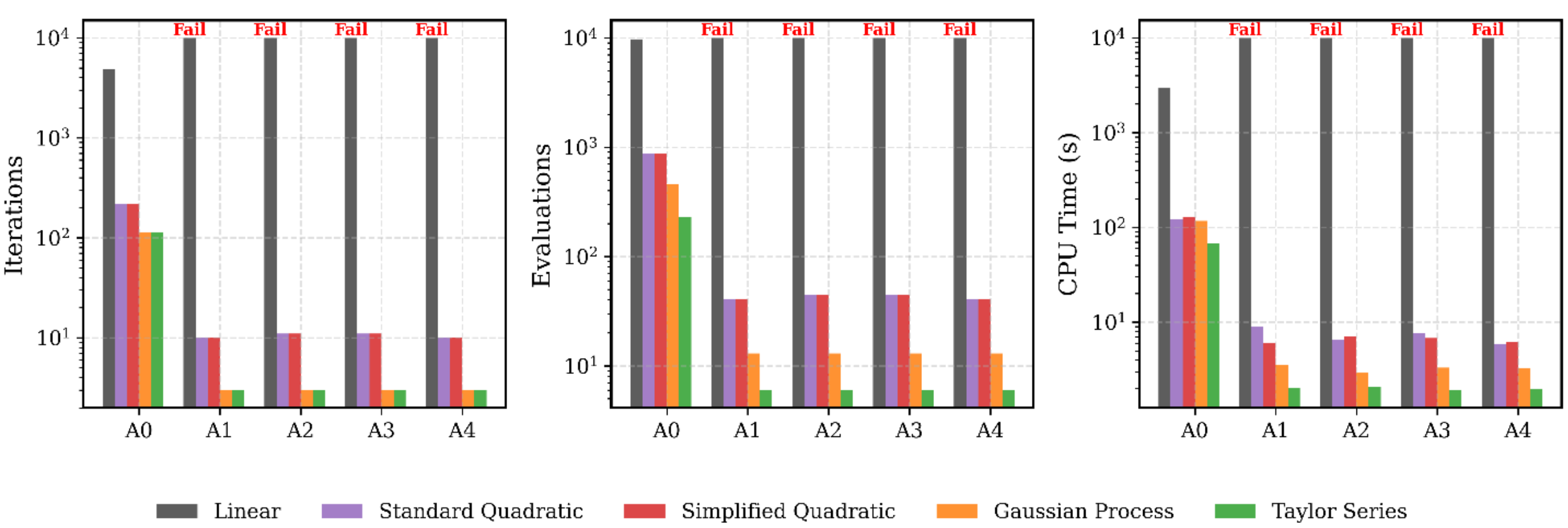

***Figure 13.*** *Comparative performance of TRF variants (A0-A4) on the Alkylation process optimisation problem. Each bar corresponds to a surrogate model. The bar with "Fail" indicated at the top shows that the specific surrogate (linear when paired with A1-A4) couldn't converge. The panels show the number of iterations (left), external black-box evaluations (middle), and computational time (right) on the y-axis (logarithmic scale). Across surrogates, A1-A4 algorithms converged in fewer iterations compared to the A0 method, demonstrating their robustness for constrained flowsheet optimisation.*

## Spring design

The tension/compression spring design problem $(n_w, n_y, n_z) = (3,1,0)$ seeks to minimise spring weight subject to nonlinear stress and deflection constraints [72]. Hessian eigenvalues showed flat curvature, reducing the dimensionality of the problem. Variants A1-A4 consistently outperformed the baseline A0 across all performance metrics (Figure 14). High-fidelity surrogates (TS and GP) proved particularly effective, solving the spring design problem in only three iterations, whereas the simplified quadratic surrogate required up to 97 iterations even with Hessian-based variants.

Across the engineering case studies, consistent patterns emerged. Hessian-informed TRF variants (A1-A4) achieved faster convergence, greater robustness, and reduced black-box evaluations relative to the baseline A0. High-fidelity surrogates (GP and TS) consistently outperformed lower-fidelity polynomial surrogates, with the combination of second-order information and accurate surrogate modelling proving critical in managing stiffness, poor conditioning, and multimodality. These findings reinforce the general advantages of Hessian projections and high-fidelity surrogates across both benchmark and applied optimisation problems.

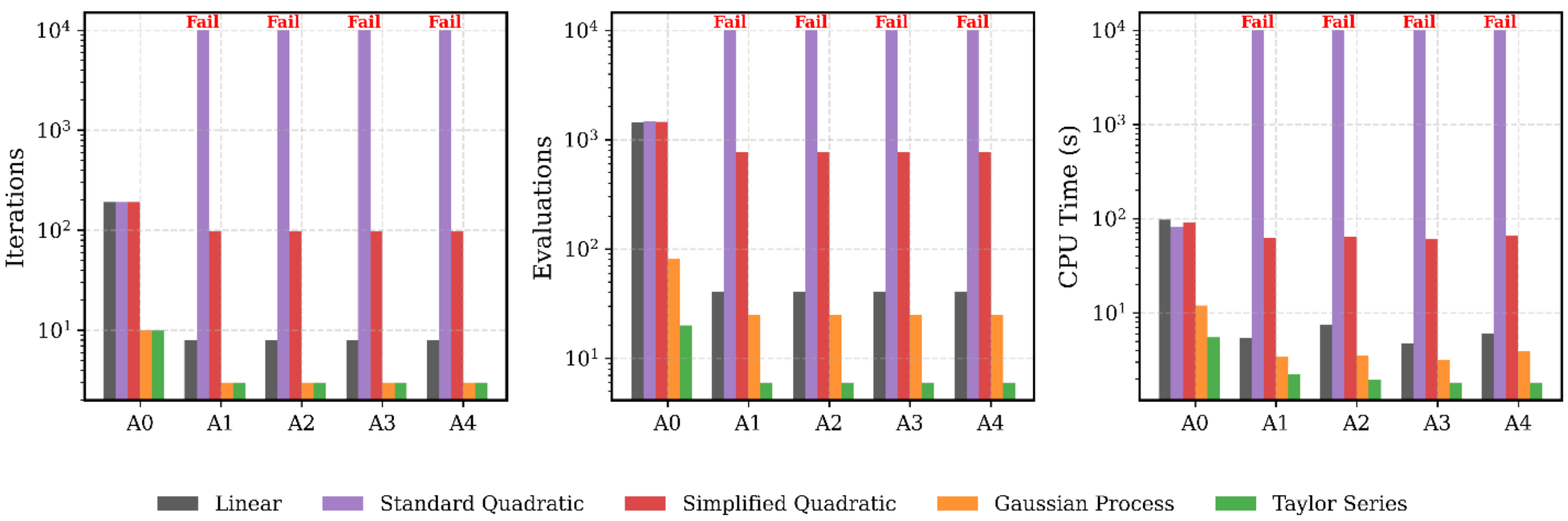

***Figure 14.*** *Comparative performance of TRF variants (A0-A4) on the spring design optimisation problem. Each bar corresponds to a surrogate model. The bar with "Fail" indicated at the top shows that the specific surrogate (standard*

*quadratic when paired with A1-A4) couldn't converge. The panels show the number of iterations (left), external black-box evaluations (middle), and computational time (right) on the y-axis (logarithmic scale). Projected-Hessian methods achieve faster convergence, demonstrating the applicability of the proposed variants across engineering design problems.*

The engineering problems displayed performance trends consistent with those observed on the broader problem set. Grey-box formulations typically converged to local minima matching those obtained by established NLP solvers (e.g., IPOPT, CONOPT). Notably, in several cases, A1-A4 reached global minima that traditional solvers failed to identify. This improvement may be influenced by the TRF framework's use of projected Hessians, which eliminates negative curvature directions and prevents premature termination at saddle points, enabling continued progress toward better solutions when within reach of the trust-region search.

While polynomial surrogates sometimes struggled in combination with Hessian-based variants, likely due to limited adaptability and lack of tuning, high-fidelity surrogates (GP, TS) consistently delivered superior performance. The enhanced robustness and reduced iteration counts observed with A1-A4 highlight the importance of exploiting curvature information. In practical terms, these improvements are highly significant for large-scale and highly nonlinear subproblems or high-cost black-box evaluations (e.g., computational fluid dynamics or detailed process simulations), where each iteration or external evaluation may require hours or days. Even modest reductions in iterations yield substantial runtime savings; reducing 100 iterations to 10, for instance, represents an order-of-magnitude efficiency gain. More broadly, our results highlight the central role of Hessian projections and high-fidelity surrogates in reducing computational effort and user intervention, thereby advancing the practical viability of grey-box optimisation.

# 6. Conclusions and Recommendations

This study improves the computational performance and robustness of the classical TRF algorithm [10] by integrating second-order information and high-fidelity surrogate forms. The modifications are validated on benchmark problems and representative engineering case studies. Among the surrogates evaluated, the TS surrogate consistently achieved the lowest number of external evaluations, proving most efficient when derivative information was available, solving all test problems under A1-A4. GP offered a strong derivative-free alternative, requiring more evaluations than TS but outperforming polynomial surrogates. Linear and quadratic surrogates solved only 72-84% of cases and required substantially more evaluations. Among the tested algorithms, the Hessian-based A2 and A3 variants were the most adaptable and robust across surrogate types (80-100% success). A4 was the fastest but less reliable with quadratic surrogates. A1 showed weaker performance than the spectrally projected Hessian methods (A2-A4) when paired with polynomial surrogates. The hybrid surrogate was least effective, solving under 20% of problems. The results also demonstrate that TRF-Hessian variants (A1-A4) not only retained robustness in reproducing local minima but also enhanced global exploration, occasionally reaching superior solutions missed by local solvers. Equally important, the number of user-tuned TRF parameters drops from five ($\epsilon_\theta$, $\epsilon_\chi$, $\epsilon_\Delta$, $\Delta^{(0)}$ and $\Delta_{min}$) in A0 to two ($\epsilon_\theta$ and $\epsilon_r$), because projected Hessians automatically control steps. We also demonstrated that the proposed grey-box TRF solvers outperform classical DFO solvers when the problem structure is partially known.

Finally, we note several limitations and avenues for further study. Projected Hessians reduce the dimensionality of the trust-region constraint; since steps are penalised differently depending on curvature sensitivity, there remains potential for multiple trust-regions or sampling-regions (for multiple heterogeneous black-box process units within a single process flowsheet) converging toward a unified optimum. The scalability of GP surrogates is a concern: in high-dimensional spaces or with many data points, GP regression can become computationally burdensome (both in fitting and evaluation). The TS surrogate approach similarly assumes smoothness and derivative availability at the current iterate; if the black-box is noisy or non-differentiable, TS-based steps may mislead the optimiser. Incorporating alternative surrogate models, such as gradient-enhanced GPs or hybridised neural-GP frameworks, may improve efficiency and accuracy. Future hybrid surrogate modelling could adaptively switch between TS and GP surrogates. Finally, our experiments focused on medium-scale problems; more work is needed to assess performance on truly large-scale, high-dimensional, or highly constrained grey-box systems. In summary, while the proposed surrogate-Hessian methods show clear promise, further research is required to improve their adaptability, reliability, and scalability for the most demanding applications.

**Data Availability Statement**

All data supporting the findings of this study are available in the associated GitHub repository [66] for reproducibility and transparency. The repository contains the model formulations/codes of all benchmark engineering case studies and the problem set, Pyomo implementations of the TRF algorithmic variants (A0-A4) and pure DFO solvers for benchmarking, and instructions for generating the numerical results reported in this work. The simplified TRF (A5) method can be accessed from the Pyomo GitHub repository [74].

Additional supplementary material documentation (Annexe A) provides complete mathematical formulations of the engineering benchmark problems (P1-P5).

## Appendix A: Global Convergence

Eason and Biegler established that the sequence of iterates generated by the classical TRF algorithm either converges to a first-order KKT point of the grey-box problem (2), provided a constraint qualification holds, or terminates unsuccessfully at a stationary point during the restoration phase after the computational budget is exhausted [21]. Note that the proof assumes $\kappa - fully\ linear$ property, regularity of the limit points and the standard NLP assumptions.

Moreover, nonlinear TRSP should always return a solution that sufficiently improves the current feasible solution.

For completeness, we summarise the convergence of the TRF algorithm to feasibility and to a first-order critical point.

**Convergence to Feasibility**

Convergence to feasibility is independent of the specific trust-region update mechanism and relies primarily on two key elements of the TRF framework: (i) the filter acceptance condition, which enforces non-dominance in the objective-infeasibility measure space, and (ii) the boundedness of the objective function $f$. The following theorem shows that, under these conditions, the algorithm generates a sequence of iterates approaching feasibility. This result assumes that the inner restoration phase terminates successfully and that the algorithm produces an infinite sequence of iterates.

**Theorem 1.** *Consider sequences* $\{\theta^{(k)}\}$ *and* $\{f^{(k)}\}$ *such that* $\theta^{(k)} \geq 0$ *and* $f^{(k)}$ *is monotonically decreasing and bounded below. Let constants* $\gamma_\theta$ *and* $\gamma_f$ *satisfy* $0 < \gamma_f < (1-\gamma_\theta) < 1$. *If, for all* $k$*,*

$$either \quad \theta^{(k+1)} \leq (1-\gamma_\theta)\theta^{(k)} \quad or \quad f^{(k)} - f^{(k+1)} \geq \gamma_f \theta^{(k)}, \qquad (A1)$$

*then* $\theta^{(k)} \to 0$ *for* $k \to \infty$.

**Proof.** Two cases follow directly from (A1):

1. If $\theta^{(k+1)} \leq (1-\gamma_\theta)\theta^{(k)}$ for all $k$, then $\theta^{(k)} \to 0$ *for* $k \to \infty$.
2. Otherwise, there exists an infinite subsequence $k'$ such that $f^{(k)} - f^{(k+1)} \geq \gamma_f \theta^{(k)}$ for all $k \geq k'$. Summing from $k = k'$ to $k = k' + N$ gives

$$\sum_{k=k'}^{k'+N} (f^{(k)} - f^{(k+1)}) \geq \gamma_f \sum_{k=k'}^{k'+N} \theta^{(k)}.$$

Thus,

$$f^{(k')} - f^{(k'+N+1)} \geq \gamma_f \sum_{k=k'}^{k'+N} \theta^{(k)}.$$

Since $f$ is bounded below and monotonically decreasing, $\sum_{k=k'}^{k'+N} \theta^{(k)} < \infty$ for $N \to \infty$, implying $\theta^{(k)} \to 0$ as $k \to \infty$.

Therefore,

$$\lim_{k\to\infty} \theta^{(k)} = 0. \qquad (A2)$$

When the filter condition (A1) holds, each accepted step is either $\theta - type$ (reducing infeasibility) or $f - type$ (reducing objective value). Hence, convergence to feasibility follows directly from the filter acceptance rule.

**Convergence to First-Order Critical Point**

Under mild assumptions [1], there exists a subsequence $\{k_{(j)}\}$ for which $x^{(k_{(j)})} \to x^*$, a first-order KKT point of the original grey-box problem (2). The TRF convergence analysis establishes that the $\theta^{(k_{(j)})}$, $\chi^{(k_{(j)})}$ and $\sigma^{(k_{(j)})}$ all approach to zero as $j \to \infty$ [10]. Consequently, $x^*$ satisfies the first-order KKT conditions, ensuring that the limit point of the feasible sequence is a locally optimal solution to the grey-box optimisation problem (2).

*Supplementary Material – Annex A*

**Himmelblau Problem (P1)**

We propose a modified Himmelblau problem that can be formulated as a constrained grey-box optimisation problem (P1) with number of variables $(n_w, n_y, n_z) = (3,2,5)$, as follows

$$\min_{w,y,z} \quad 5.3578547y_1 + 0.8356891z_1w_5 + 37.2932239z_1 - 40792.141, \tag{P1}$$

$$
\begin{aligned}
s.t. \quad & h_1(w,y,z) := z_6 = 85.334407 + 0.0056858y_2 + 0.00026z_1z_4 - 0.0022053w_3w_5, \\
& h_2(w,y,z) := z_7 \\
& \qquad = 80.51249 + 0.0071317y_2 + 0.0029955z_1w_2 - 0.0021813{w_3}^2, \\
& h_3(w,y,z) := z_8 \\
& \qquad = 9.300961 + 0.0047026w_3w_5 + 0.0012547z_1w_3 \\
& \qquad - 0.0019085w_3z_4, \\
& y_1 = d_1(w) := w_3^2, \\
& y_2 = d_2(w) := w_2w_5, \\
& 78 \leq z_1 \leq 102, \\
& 33 \leq w_2 \leq 45, \\
& 27 \leq w_i \leq 45, \forall i \in \{3,5\}, \\
& 27 \leq z_4 \leq 45, \\
& 0 \leq z_6 \leq 92, \\
& 90 \leq z_7 \leq 110, \\
& 20 \leq z_8 \leq 25.
\end{aligned}
$$

**Liquid-liquid Extraction Column (P2)**

We consider a the optimisation of flowrates in liquid-liquid extraction column as a constrained grey-box optimisation problem (P2) with variable dimensions $(n_w, n_y, n_z) = (2,1,2)$:

$$\min_{w,y,z} \quad -w_2z_1, \tag{P2}$$

$$
\begin{aligned}
s.t. \quad & g_1(w,y,z) := w_1 + w_2 \leq 0.2, \\
& h_1(w,y,z) := z_2 - (1.5w_1/w_2) = 0, \\
& h_2(w,y,z) := z_1\left(1 - z_2e^{y_1(1-z_2)}\right) - z_2\left(1 - e^{y_1(1-z_2)}\right) = 0, \\
& y_1 = d_1(w) := 4.81(w_1/w_2)^{0.24}, \\
& 0.05 \leq w_1 \leq 0.25, \\
& 0.05 \leq w_2 \leq 0.3,
\end{aligned}
$$

$0 \leq z_1 \leq 10,$

$0 \leq z_2 \leq 10.$

**Pressure Vessel Design (P3)**

We consider a modified formulation of the pressure vessel design problem as a constrained grey-box optimisation problem (P3) with number of variables $(n_w, n_y, n_z) = (4,1,0)$, as follows

$$\min_{w,y,z} \quad y_1, \tag{P3}$$

$$
\begin{aligned}
s.t. \quad & g_1(w,y,z) := -w_1 + 0.0193w_3 \leq 0,\\
& g_2(w,y,z) := -w_2 + 0.0095w_3 \leq 0,\\
& g_3(w,y,z) := -\pi {w_3}^2 w_4 - \frac{4}{3}\pi {w_3}^3 + 1296000 \leq 0,\\
& y_1 = d_1(w) := 0.6224(w_1 w_3 w_4) + 1.7781 w_2 {w_3}^2 + 3.1661 {w_1}^2 w_4 + 19.84 {w_1}^2 w_3,\\
& 0.0625 \leq \ w_1 \leq 6.1875,\\
& 0.0625 \leq w_2 \leq 6.1875,\\
& 10 \leq w_3 \leq 200,\\
& 10 \leq w_4 \leq 200.
\end{aligned}
$$

**Alkylation Process (P4)**

We reformulate Alkylation process problem into a constrained grey-box optimisation problem (P4) with number of variables $(n_w, n_y, n_z) = (1,1,9)$, as follows

$$\max_{w,y,z} \quad 0.063z_4 z_7 - 5.04z_1 - 0.035z_2 - 10z_3 - 3.36z_5, \tag{P4}$$

$$
\begin{aligned}
s.t. \quad & h_1(w,y,z) := z_4 - z_1 y_1 = 0,\\
& h_2(w,y,z) := z_7 - 86.35 - 1.098w_1 + 0.038{w_1}^2 - 0.325(z_6 - 89) = 0,\\
& h_3(w,y,z) := z_8 - 35.82 + 0.222z_9 = 0,\\
& h_4(w,y,z) := z_9 + 133 - 3z_7 = 0,\\
& h_5(w,y,z) := w_1 z_1 - z_2 - z_5 = 0,\\
& h_6(w,y,z) := z_5 - 1.22z_4 + z_1 = 0,\\
& h_7(w,y,z) := z_6(z_4 z_8 + 1000z_3) - 98000z_3 = 0,\\
& y_1 = d_1(w) := 1.12 + 1.2167w_1 + 0.0067{w_1}^2,\\
& 0 \leq z_1 \leq 2000,\\
& 0 \leq z_2 \leq 16000,
\end{aligned}
$$

$0 \leq z_3 \leq 120,$

$0 \leq z_4 \leq 5000,$

$0 \leq z_5 \leq 2000,$

$85 \leq z_6 \leq 93,$

$90 \leq z_7 \leq 95,$

$1.2 \leq z_8 \leq 4,$

$145 \leq z_9 \leq 162,$

$3 \leq w_1 \leq 12.$

**Spring Design (P5)**

We reformulate spring design problem into a constrained grey-box optimisation problem (P5) with number of variables $(n_w, n_y, n_z) = (3,1,0)$, as follows

$$\min_{w,y,z} \quad y_1, \tag{P5}$$

$$\text{s.t.} \quad g_1(w,y,z) := 1 - \frac{{w_2}^3 w_3}{71785 {w_1}^4} \leq 0,$$

$$g_2(w,y,z) := \frac{4{w_2}^2 - w_1 w_2}{12566(w_2 {w_1}^3 - {w_1}^4)} + \frac{1}{5108 {w_1}^2} - 1 \leq 0,$$

$$g_3(w,y,z) := 1 - \frac{140.45 w_1}{{w_2}^2 w_3} \leq 0,$$

$$y_1 = d_1(w) := (w_3 + 2) w_2 {w_1}^2,$$

$$0.05 \leq w_1 \leq 2,$$

$$0.25 \leq w_2 \leq 1.3,$$

$$2 \leq w_3 \leq 15.$$